\documentclass[12pt,thmsa]{article}%
\usepackage{sw20bams}
\usepackage{amsmath}
\usepackage{amsfonts}
\usepackage{amssymb}
\usepackage{graphicx}%
\providecommand{\U}[1]{\protect\rule{.1in}{.1in}}
\begin{document}

\author{Steven R. Finch and Antonia J. Jones}
\title{Random Spherical Triangles}
\date{December 21, 2015}
\maketitle

\begin{abstract}
Let $\Delta$ be a random spherical triangle (meaning that vertices are
independent and uniform on the unit sphere). A closed-form expression for the
area density of $\Delta$ has been known since 1867; a complicated integral
expression for the perimeter density was found in 1994. Does there exist a
closed-form expression for the latter? We attempt to answer this question from
several directions. An outcome of our work is the exact value of the perimeter
density at the point $\pi$.

\end{abstract}

\footnotetext{Copyright \copyright \ 2010, 2015 by Steven R. Finch. All rights
reserved.}A\ spherical triangle $\Delta$ is a region enclosed by three great
circles on the unit sphere; a great circle is a circle whose center is at the
origin. The sides of $\Delta$ are arcs of great circles and have length $a$,
$b$, $c$. Each of these is $\leq\pi$. The angle $\alpha$ opposite side $a$ is
the dihedral angle between the two planes passing through the origin and
determined by arcs $b$, $c$. The angles $\beta$, $\gamma$ opposite sides $b$,
$c$ are similarly defined. Each of these is $\leq\pi$ too.

Define a \textbf{primal triangle} to be a random spherical triangle, obtained
by selecting three independent uniformly distributed points $A$, $B$, $C$ on
the sphere to be vertices. Define a \textbf{dual triangle} to be a random
spherical triangle, obtained by selecting three independent uniformly
distributed great circles on the sphere to be sides. More precisely, starting
with independent uniform points $A^{\prime}$, $B^{\prime}$, $C^{\prime}$ on
the sphere, a dual triangle has vertices
\[%
\begin{array}
[c]{ccccc}%
A=\dfrac{B^{\prime}\times C^{\prime}}{\left\|  B^{\prime}\times C^{\prime
}\right\|  }, &  & B=\dfrac{A^{\prime}\times C^{\prime}}{\left\|  A^{\prime
}\times C^{\prime}\right\|  }, &  & C=\dfrac{A^{\prime}\times B^{\prime}%
}{\left\|  A^{\prime}\times B^{\prime}\right\|  }.
\end{array}
\]
Hence, while its vertices are not independent, the poles of a dual triangle are.

Let $\Delta$ be a primal triangle. The area $\sigma=\alpha+\beta+\gamma-\pi$
of $\Delta$ satisfies $0\leq\sigma\leq2\pi$. The perimeter $\tau=a+b+c$ of
$\Delta$ satisfies $0\leq\tau\leq2\pi$. Expressions for the trivariate density
of $(\alpha,\beta,\gamma)$ and the trivariate density of $(a,b,c)$ are known
\cite{Ms} but do not give useful insight into the distributions of $\sigma$
and $\tau$. Crofton \& Exhumatus \cite{CE} determined the density for $\sigma
$:
\[
-\frac{(x^{2}-4\pi x+3\pi^{2}-6)\cos(x)-6(x-2\pi)\sin(x)-2(x^{2}-4\pi
x+3\pi^{2}+3)}{16\pi\cos(x/2)^{4}}
\]
for $0<x<2\pi$; this formula remained obscure until it was cited in a recent
paper \cite{HDB}. Unpublished work of J. N. Boots (mentioned in \cite{Ms}) was
also apparently relevant. Jones \&\ Benyon-Tinker \cite{ABT1, ABT2} determined
the density for $\tau$:
\[
\frac1{4\pi}%
{\displaystyle\int\limits_{0}^{x/2}}
\frac{E\left(  \sin\left(  \frac t2\right)  \right)  -\cos\left(  \frac
{x-t}2\right)  ^{2}K\left(  \sin\left(  \frac t2\right)  \right)  }{\sqrt
{\cos\left(  \frac t2\right)  ^{2}-\cos\left(  \frac{x-t}2\right)  ^{2}}}%
\sin(t)\,dt
\]
for $0<x<2\pi$, where
\[%
\begin{array}
[c]{c}%
K(\zeta)=%
{\displaystyle\int\limits_{0}^{\pi/2}}
\dfrac1{\sqrt{1-\zeta^{2}\sin(\theta)^{2}}}\,d\theta=%
{\displaystyle\int\limits_{0}^{1}}
\dfrac1{\sqrt{(1-t^{2})(1-\zeta^{2}t^{2})}}\,dt,\\
E(\zeta)=%
{\displaystyle\int\limits_{0}^{\pi/2}}
\sqrt{1-\zeta^{2}\sin(\theta)^{2}}\,d\theta=%
{\displaystyle\int\limits_{0}^{1}}
\sqrt{\dfrac{1-\zeta^{2}t^{2}}{1-t^{2}}}\,dt
\end{array}
\]
are complete elliptic integrals of the first and second kind. We wonder: does
there exist a closed-form expression for this latter density? A direct
evaluation of the integral does not seem possible, yet conceivably a different
geometric argument might yield a more accessible formula. We attempt to answer
this question from several directions. Also, it is clear that the density is
zero at $x=0$ and diverges to infinity at $x=2\pi$. Numerically the density
$\approx3\sqrt{2}/32$ to high precision at $x=\pi$, but a proof via the
preceding is not known. An outcome of our work is a new formula that gives the
exact value as predicted.

\section{Two Coordinate Systems}

We define two coordinate systems on the unit sphere that will help in our
study of triangular area and density.

\subsection{\label{PrimalCoord}Primal Coordinates}

Without loss of generality, let $A=(1,0,0)$ and $B=(\cos(\kappa),\sin
(\kappa),0)$ in $xyz$ coordinates. We wish to locate the unique point $C$ in
the upper hemisphere so that the triangle $ABC$ satisfies $\alpha=\theta$,
$b=\rho$, $c=\kappa$. See Figure 1.
\begin{figure}[ptb]%
\centering
\includegraphics[
height=3.6097in,
width=3.5898in
]%
{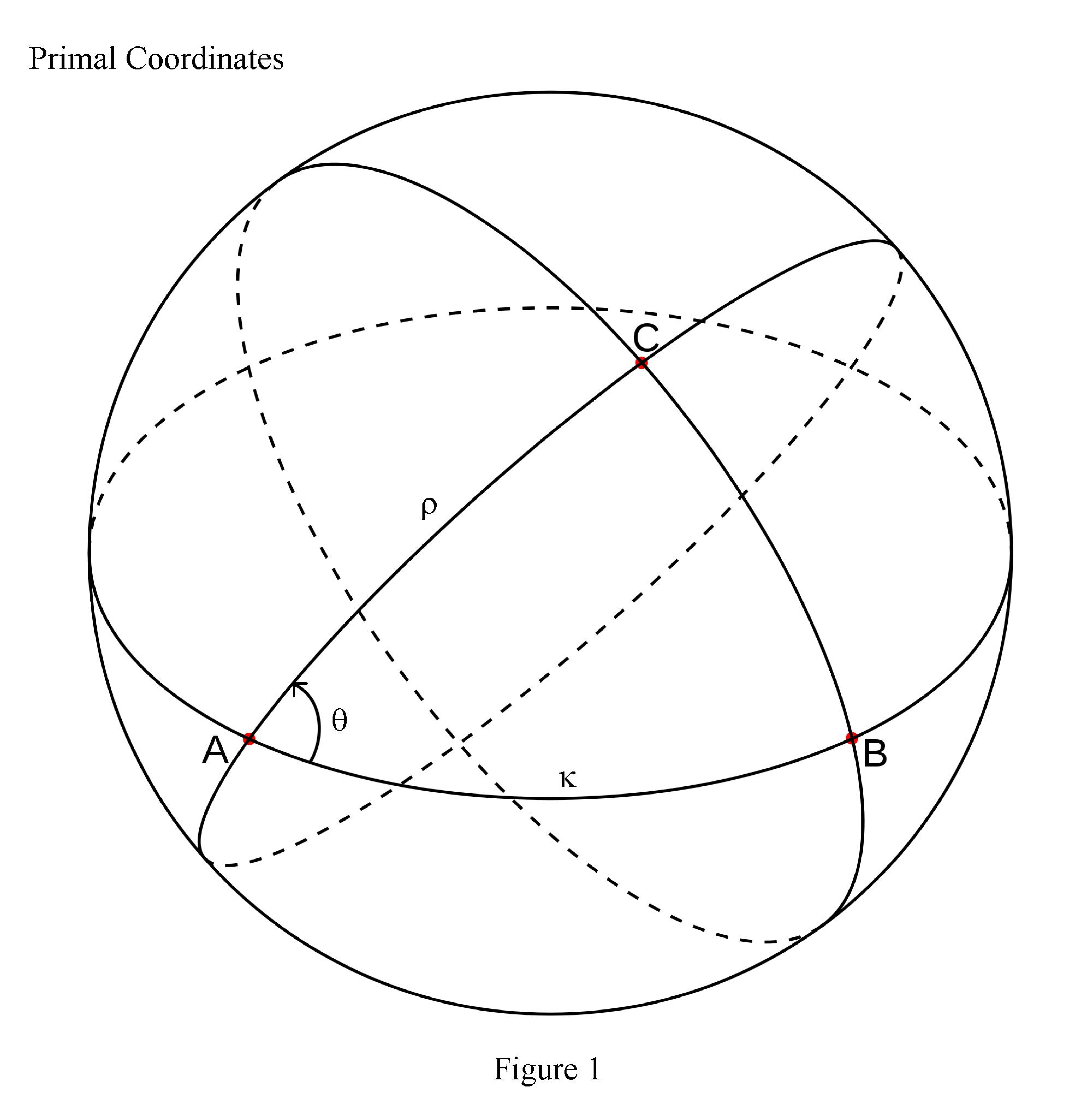}%
\end{figure}
The parameters $\rho$, $\theta$ are regarded as varying while the parameter
$\kappa$ is fixed. Think of rotating the equatorial disk in space so that the
vector $(1,0,0)$ remains fixed and the vector $(0,1,0)$ moves toward $(0,0,1)$
through the angle $\theta$. The rotation matrix performing this motion is
\cite{Wi}
\[
R=\left(
\begin{array}
[c]{ccc}%
1 & 0 & 0\\
0 & \cos(\theta) & -\sin(\theta)\\
0 & \sin(\theta) & \cos(\theta)
\end{array}
\right)
\]
and
\[
R\left(
\begin{array}
[c]{c}%
\cos(\rho)\\
\sin(\rho)\\
0
\end{array}
\right)  =\left(
\begin{array}
[c]{c}%
\cos(\rho)\\
\sin(\rho)\cos(\theta)\\
\sin(\rho)\sin(\theta)
\end{array}
\right)  ,
\]
which gives the point $C$. The three-dimensional transformation
\[
\left(
\begin{array}
[c]{c}%
r\\
\rho\\
\theta
\end{array}
\right)  \mapsto\left(
\begin{array}
[c]{c}%
x\\
y\\
z
\end{array}
\right)  =\left(
\begin{array}
[c]{c}%
r\cos(\rho)\\
r\sin(\rho)\cos(\theta)\\
r\sin(\rho)\sin(\theta)
\end{array}
\right)
\]
has Jacobian determinant
\[
\left\vert \frac{\partial(x,y,z)}{\partial(r,\rho,\theta)}\right\vert
=\left\vert
\begin{array}
[c]{ccc}%
\cos(\rho) & -r\sin(\rho) & 0\\
\sin(\rho)\cos(\theta) & r\cos(\rho)\cos(\theta) & -r\sin(\rho)\sin(\theta)\\
\sin(\rho)\sin(\theta) & r\cos(\rho)\sin(\theta) & r\sin(\rho)\cos(\theta)
\end{array}
\right\vert =r^{2}\sin(\rho)
\]
which implies that the area element in primal $\rho\theta$ coordinates is
$\sin(\rho)d\rho\,d\theta$.

\subsection{\label{DualCoord}Dual Coordinates}

Without loss of generality, let $A=(1,0,0)$ and $B=(\cos(\rho),\sin(\rho),0)$
in $xyz$ coordinates. It seems (at first glance) that we should locate the
unique point $C$ in the upper hemisphere so that the triangle $ABC$ satisfies
$\alpha=\kappa$, $\beta=\theta$, $c=\rho$. See Figure 2.
\begin{figure}[ptb]%
\centering
\includegraphics[
height=3.5457in,
width=3.5267in
]%
{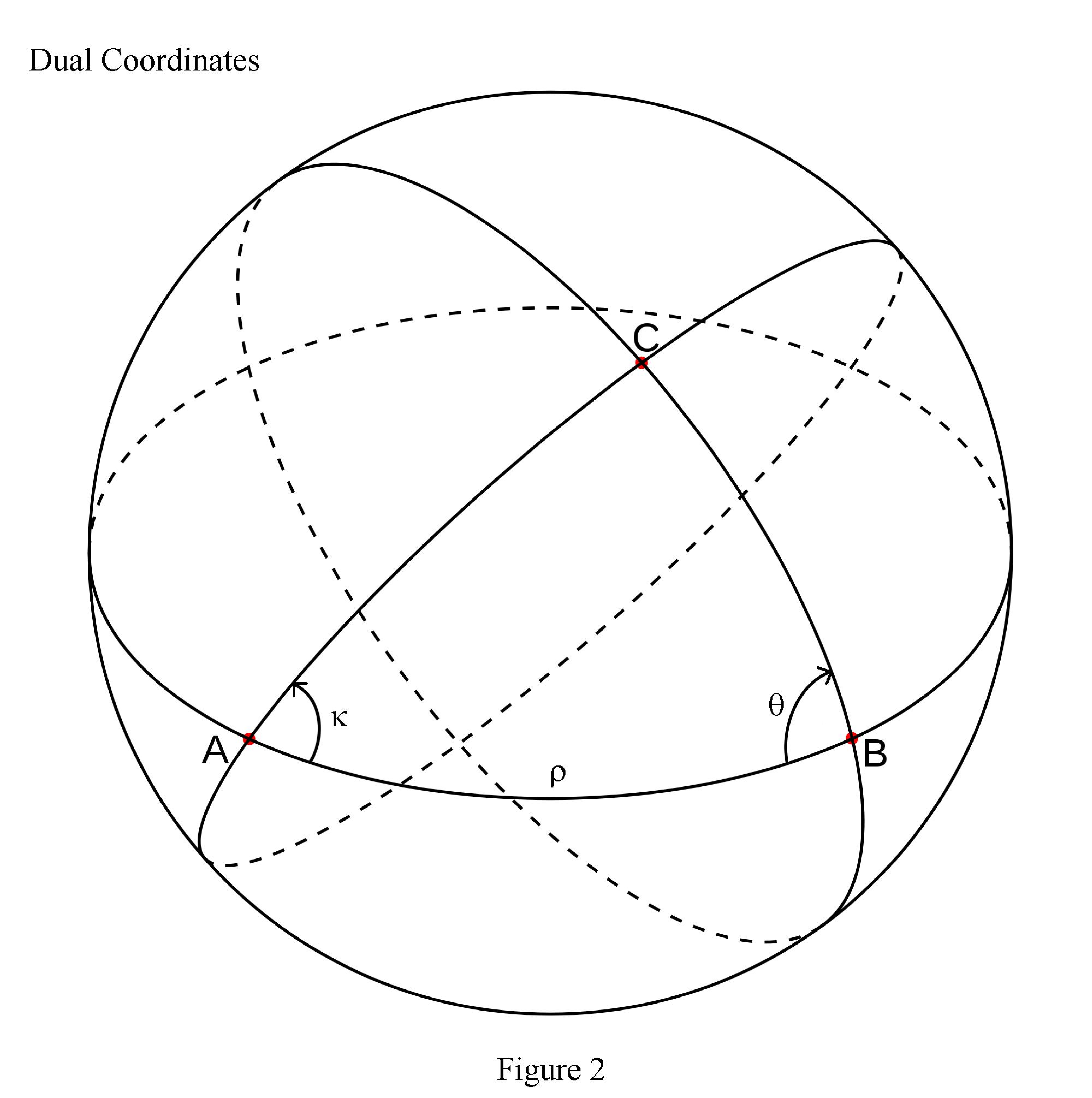}%
\end{figure}
The parameters $\rho$, $\theta$ are regarded as varying while the parameter
$\kappa$ is fixed.

Let us examine the great circle containing $A$, $C$. It must also contain the
point $(0,\cos(\kappa),\sin(\kappa))$, since this is the image of $(0,1,0)$
after rotation through angle $\kappa$. Hence a normal vector is
$V=(1,0,0)\times(0,\cos(\kappa),\sin(\kappa))=(0,-\sin(\kappa),\cos(\kappa))$.

Let us examine the great circle containing $B$, $C$. Think of rotating the
equatorial disk in space so that the vector $(\cos(\rho),\sin(\rho),0)$
remains fixed and the vector $(\sin(\rho),-\cos(\rho),0)$ moves toward
$(0,0,1)$ through the angle $\theta$. The rotation matrix performing this
motion is \cite{Wi}
\[
S=\left(
\begin{array}
[c]{ccc}%
\cos(\rho)^{2}+(1-\cos(\rho)^{2})\cos(\theta) & \cos(\rho)\sin(\rho
)(1-\cos(\theta)) & -\sin(\rho)\sin(\theta)\\
\cos(\rho)\sin(\rho)(1-\cos(\theta)) & \sin(\rho)^{2}+(1-\sin(\rho)^{2}%
)\cos(\theta) & \cos(\rho)\sin(\theta)\\
\sin(\rho)\sin(\theta) & -\cos(\rho)\sin(\theta) & \cos(\theta)
\end{array}
\right)
\]
and
\[
S\left(
\begin{array}
[c]{c}%
\sin(\rho)\\
-\cos(\rho)\\
0
\end{array}
\right)  =\left(
\begin{array}
[c]{c}%
\sin(\rho)\cos(\theta)\\
-\cos(\rho)\cos(\theta)\\
\sin(\theta)
\end{array}
\right)  .
\]
For example, if $\rho=\pi/2$, the image of $(1,0,0)$ after rotation through
angle $\theta$ is $(\cos(\theta),0,\sin(\theta))$. As another example, if
$\rho=0$, the image of $(0,-1,0)$ after rotation through angle $\theta$ is
$(0,-\cos(\theta),\sin(\theta))$. Hence the great circle must contain the
point $(\sin(\rho)\cos(\theta),-\cos(\rho)\cos(\theta),\sin(\theta))$ and a
normal vector is
\[
W=\left(
\begin{array}
[c]{c}%
\cos(\rho)\\
\sin(\rho)\\
0
\end{array}
\right)  \times\left(
\begin{array}
[c]{c}%
\sin(\rho)\cos(\theta)\\
-\cos(\rho)\cos(\theta)\\
\sin(\theta)
\end{array}
\right)  =\left(
\begin{array}
[c]{c}%
\sin(\rho)\sin(\theta)\\
-\cos(\rho)\sin(\theta)\\
-\cos(\theta)
\end{array}
\right)  .
\]
The point $C$ is orthogonal to the two normal vectors and at unit distance
from the origin, equivalently, $C=(V\times W)/\left\|  V\times W\right\|  $.

Now, in fact, this is more than what is required. We need only (on second
glance) specify the great circle containing $B$, $C$ and this is done via
locating $W$ or $-W$. The three-dimensional transformation
\[
\left(
\begin{array}
[c]{c}%
r\\
\rho\\
\theta
\end{array}
\right)  \mapsto\left(
\begin{array}
[c]{c}%
x\\
y\\
z
\end{array}
\right)  =\left(
\begin{array}
[c]{c}%
-r\sin(\rho)\sin(\theta)\\
r\cos(\rho)\sin(\theta)\\
r\cos(\theta)
\end{array}
\right)
\]
has Jacobian determinant $r^{2}\sin(\theta)$, which implies that the area
element in dual $\rho\theta$ coordinates is $\sin(\theta)d\rho\,d\theta$.
Perhaps this is obvious by duality. It is good, however, to see the supporting
geometric details.

\section{Four Approaches}

We illustrate using four different trigonometric identities and the above two
coordinate systems. More possible approaches will be mentioned in a later section.

\subsection{\label{PrimalArea}Primal Area}

As in section [\ref{PrimalCoord}], assume that the triangle $ABC$ satisfies
$\alpha=\theta$, $b=\rho$, $c=\kappa$. These three parameters are related to
area $\sigma$ as follows:
\[
\tan\left(  \frac\rho2\right)  =\cot\left(  \frac\kappa2\right)  \csc\left(
\theta-\frac\sigma2\right)  \sin\left(  \frac\sigma2\right)  .
\]
A\ proof appears in section [\ref{TrigIdentPrimal1}]. For fixed $\sigma$ and
$\kappa$, define
\[
f(\theta)=\left\{
\begin{array}
[c]{lll}%
\pi &  & \text{if }0\leq\theta<\sigma/2,\\
2\arctan\left[  \cot\left(  \dfrac\kappa2\right)  \csc\left(  \theta
-\dfrac\sigma2\right)  \sin\left(  \dfrac\sigma2\right)  \right]  &  &
\text{if }\sigma/2\leq\theta\leq\pi
\end{array}
\right.
\]
then the conditional probability, given $c$, is
\begin{align*}
\operatorname*{P}\left\{  \text{area}\leq\sigma\left|  c=\kappa\right.
\right\}   &  =\frac1{2\pi}%
{\displaystyle\int\limits_{0}^{\pi}}
{\displaystyle\int\limits_{0}^{f(\theta)}}
\sin(\rho)d\rho\,d\theta\\
\  &  =\frac1{2\pi}%
{\displaystyle\int\limits_{0}^{\sigma/2}}
{\displaystyle\int\limits_{0}^{\pi}}
\sin(\rho)d\rho\,d\theta+\frac1{2\pi}%
{\displaystyle\int\limits_{\sigma/2}^{\pi}}
{\displaystyle\int\limits_{0}^{f(\theta)}}
\sin(\rho)d\rho\,d\theta\\
\  &  =\frac1{2\pi}\left(  \sigma+%
{\displaystyle\int\limits_{\sigma/2}^{\pi}}
(1-\cos(f(\theta)))d\theta\right)  .
\end{align*}
This result can be experimentally verified by generating many primal triangles
$ABC$ with $c=\kappa$, and then plotting all pairs $(\theta,\rho)$
corresponding to triangles with area $\leq\sigma$. The scatterplot fills the
region $[0,\pi]\times[0,\pi]$ except for the portion lying above the curve
$\rho=f(\theta)$.

We will later discuss [\ref{Success}] how the unconditional probability
$\operatorname*{P}\left\{  \text{area}\leq\sigma\right\}  $ is evaluated
exactly, for arbitrary $\sigma$. The method is quite long and intricate.

Here is a quick method for evaluating not probability, but instead density, at
$\sigma=\pi$. We start with the conditional density
\begin{align*}
\frac d{d\sigma}\frac1{2\pi}\left(  \sigma+%
{\displaystyle\int\limits_{\sigma/2}^{\pi}}
(1-\cos(f(\theta)))d\theta\right)   &  =\frac1{2\pi}\left(  1-\tfrac12\left[
1-\cos\left(  f\left(  \tfrac\sigma2\right)  \right)  \right]  +%
{\displaystyle\int\limits_{\sigma/2}^{\pi}}
\frac d{d\sigma}(1-\cos(f(\theta)))d\theta\right) \\
\  &  =\frac1{2\pi}%
{\displaystyle\int\limits_{\sigma/2}^{\pi}}
\sin(f(\theta))g(\theta)d\theta
\end{align*}
where
\[
\sin(f(\theta))=\frac{2\tan\left(  \tfrac\kappa2\right)  \sin\left(
\theta-\tfrac\sigma2\right)  \sin\left(  \tfrac\sigma2\right)  }{\tan\left(
\tfrac\kappa2\right)  ^{2}\sin\left(  \theta-\tfrac\sigma2\right)  ^{2}%
+\sin\left(  \tfrac\sigma2\right)  ^{2}}
\]
since
\[
\sin(2\arctan(\zeta))=\frac{2\zeta}{1+\zeta^{2}},
\]
and where
\[
g(\theta)=\frac{df}{d\sigma}=\frac{\tan\left(  \tfrac\kappa2\right)
\sin(\theta)}{\tan\left(  \tfrac\kappa2\right)  ^{2}\sin\left(  \theta
-\tfrac\sigma2\right)  ^{2}+\sin\left(  \tfrac\sigma2\right)  ^{2}}.
\]
It follows that the unconditional density is
\[
\frac1{2\pi}%
{\displaystyle\int\limits_{0}^{\pi}}
{\displaystyle\int\limits_{\sigma/2}^{\pi}}
\frac{\tan\left(  \tfrac\kappa2\right)  ^{2}\sin\left(  \theta-\tfrac
\sigma2\right)  \sin\left(  \tfrac\sigma2\right)  \sin\left(  \theta\right)
}{\left[  \tan\left(  \tfrac\kappa2\right)  ^{2}\sin\left(  \theta
-\tfrac\sigma2\right)  ^{2}+\sin\left(  \tfrac\sigma2\right)  ^{2}\right]
^{2}}\sin(\kappa)d\theta\,d\kappa
\]
because the density for $\kappa$ is $\sin(\kappa)/2.$ By the half-angle
formula for tangent, this is the same as
\[
\frac1{2\pi}%
{\displaystyle\int\limits_{0}^{\pi}}
{\displaystyle\int\limits_{\sigma/2}^{\pi}}
\frac{\left(  1-\cos(\kappa)\right)  \left(  1+\cos(\kappa)\right)
\sin\left(  \theta-\tfrac\sigma2\right)  \sin\left(  \tfrac\sigma2\right)
\sin\left(  \theta\right)  }{\left[  \left(  1-\cos(\kappa)\right)
\sin\left(  \theta-\tfrac\sigma2\right)  ^{2}+\left(  1+\cos(\kappa)\right)
\sin\left(  \tfrac\sigma2\right)  ^{2}\right]  ^{2}}\sin(\kappa)d\theta
\,d\kappa.
\]
In the special case when $\sigma=\pi$, this becomes
\begin{align*}
&  \ -\frac1{2\pi}%
{\displaystyle\int\limits_{0}^{\pi}}
{\displaystyle\int\limits_{\pi/2}^{\pi}}
\frac{\cos\left(  \theta\right)  \sin\left(  \theta\right)  \left(
1-\cos(\kappa)^{2}\right)  \sin(\kappa)}{\left[  1+\cos(\kappa)+\left(
1-\cos(\kappa)\right)  \cos\left(  \theta\right)  ^{2}\right]  ^{2}}%
d\theta\,d\kappa\\
\  &  =-\frac1{2\pi}%
{\displaystyle\int\limits_{-1}^{1}}
{\displaystyle\int\limits_{-1}^{0}}
\frac{u\left(  1-v^{2}\right)  }{\left[  1+v+\left(  1-v\right)  u^{2}\right]
^{2}}du\,dv=\frac1{4\pi}%
\end{align*}
consistent with Crofton \& Exhumatus. No analogous simplication seems to
occur, for example, when $\sigma=\pi/2$ or $\sigma=3\pi/2$.

\subsection{Dual Perimeter}

As in section [\ref{DualCoord}], assume that the triangle $ABC$ satisfies
$\alpha=\kappa$, $\beta=\theta$, $c=\rho$. These three parameters are related
to perimeter $\tau$ as follows:
\[
\tan\left(  \frac\theta2\right)  =\cot\left(  \frac\kappa2\right)  \sin\left(
\frac\tau2-\rho\right)  \csc\left(  \frac\tau2\right)  .
\]
A\ proof appears in section [\ref{TrigIdentDual1}]. For fixed $\tau$ and
$\kappa$, define
\[
f(\rho)=\left\{
\begin{array}
[c]{lll}%
2\arctan\left[  \cot\left(  \frac\kappa2\right)  \sin\left(  \frac\tau
2-\rho\right)  \csc\left(  \frac\tau2\right)  \right]  &  & \text{if }%
0\leq\rho\leq\tau/2,\\
0 &  & \text{if }\tau/2<\rho\leq\pi
\end{array}
\right.
\]
then the conditional probability, given $\alpha$, is
\begin{align*}
\operatorname*{P}\left\{  \text{perimeter}\leq\tau\left|  \alpha
=\kappa\right.  \right\}   &  =\frac1{2\pi}%
{\displaystyle\int\limits_{0}^{\pi}}
{\displaystyle\int\limits_{0}^{f(\rho)}}
\sin(\theta)d\theta\,d\rho\\
\  &  =\frac1{2\pi}%
{\displaystyle\int\limits_{0}^{\tau/2}}
{\displaystyle\int\limits_{0}^{f(\rho)}}
\sin(\theta)d\theta\,d\rho+\frac1{2\pi}%
{\displaystyle\int\limits_{\tau/2}^{\pi}}
{\displaystyle\int\limits_{0}^{0}}
\sin(\theta)d\theta\,d\rho\\
\  &  =\frac1{2\pi}%
{\displaystyle\int\limits_{0}^{\tau/2}}
(1-\cos(f(\rho)))d\rho.
\end{align*}
This result can be experimentally verified by generating many dual triangles
$ABC$ with $\alpha=\kappa$, and then plotting all pairs $(\rho,\theta)$
corresponding to triangles with perimeter $\leq\tau$. The scatterplot fills
the region $[0,\pi]\times[0,\pi]$ except for the portion lying above the curve
$\theta=f(\rho)$.

Since (dual perimeter) = ($2\pi-$ primal area), it is not surprising that
conditional probabilities are so similar.

For completeness' sake, let us compute the conditional density
\begin{align*}
\frac d{d\tau}\frac1{2\pi}%
{\displaystyle\int\limits_{0}^{\tau/2}}
(1-\cos(f(\rho)))d\rho &  =\frac1{2\pi}\left(  \tfrac12\left[  1-\cos\left(
f\left(  \tfrac\tau2\right)  \right)  \right]  +%
{\displaystyle\int\limits_{0}^{\tau/2}}
\frac d{d\tau}(1-\cos(f(\rho)))d\rho\right) \\
\  &  =\frac1{2\pi}%
{\displaystyle\int\limits_{0}^{\tau/2}}
\sin(f(\rho))g(\rho)d\rho
\end{align*}
where
\[
\sin(f(\rho))=\frac{2\tan\left(  \tfrac\kappa2\right)  \sin\left(  \tfrac
\tau2-\rho\right)  \sin\left(  \tfrac\tau2\right)  }{\tan\left(  \tfrac
\kappa2\right)  ^{2}\sin\left(  \tfrac\tau2\right)  ^{2}+\sin\left(
\tfrac\tau2-\rho\right)  ^{2}},
\]
\[
g(\rho)=\frac{df}{d\tau}=\frac{\tan\left(  \tfrac\kappa2\right)  \sin(\rho
)}{\tan\left(  \tfrac\kappa2\right)  ^{2}\sin\left(  \tfrac\tau2\right)
^{2}+\sin\left(  \tfrac\tau2-\rho\right)  ^{2}}.
\]
It follows that the unconditional density is
\begin{align*}
&  \frac1{2\pi}%
{\displaystyle\int\limits_{0}^{\pi}}
{\displaystyle\int\limits_{0}^{\tau/2}}
\frac{\tan\left(  \tfrac\kappa2\right)  ^{2}\sin\left(  \tfrac\tau
2-\rho\right)  \sin\left(  \tfrac\tau2\right)  \sin\left(  \rho\right)
}{\left[  \tan\left(  \tfrac\kappa2\right)  ^{2}\sin\left(  \tfrac
\tau2\right)  ^{2}+\sin\left(  \tfrac\tau2-\rho\right)  ^{2}\right]  ^{2}}%
\sin(\kappa)d\rho\,d\kappa\\
&  =\frac1{2\pi}%
{\displaystyle\int\limits_{0}^{\pi}}
{\displaystyle\int\limits_{0}^{\tau/2}}
\frac{\left(  1-\cos(\kappa)\right)  \left(  1+\cos(\kappa)\right)
\sin\left(  \tfrac\tau2-\rho\right)  \sin\left(  \tfrac\tau2\right)
\sin\left(  \rho\right)  }{\left[  \left(  1-\cos(\kappa)\right)  \sin\left(
\tfrac\tau2\right)  ^{2}+\left(  1+\cos(\kappa)\right)  \sin\left(  \tfrac
\tau2-\rho\right)  ^{2}\right]  ^{2}}\sin(\kappa)d\rho\,d\kappa
\end{align*}
because the density for $\kappa$ is $\sin(\kappa)/2.$ In the special case when
$\tau=\pi$, this becomes
\begin{align*}
&  \ \ \frac1{2\pi}%
{\displaystyle\int\limits_{0}^{\pi}}
{\displaystyle\int\limits_{0}^{\pi/2}}
\frac{\cos\left(  \rho\right)  \sin\left(  \rho\right)  \left(  1-\cos
(\kappa)^{2}\right)  \sin(\kappa)}{\left[  1-\cos(\kappa)+\left(
1+\cos(\kappa)\right)  \cos\left(  \rho\right)  ^{2}\right]  ^{2}}%
d\rho\,d\kappa\\
\  &  =\frac1{2\pi}%
{\displaystyle\int\limits_{-1}^{1}}
{\displaystyle\int\limits_{0}^{1}}
\frac{u\left(  1-v^{2}\right)  }{\left[  1-v+\left(  1+v\right)  u^{2}\right]
^{2}}du\,dv=\frac1{4\pi}%
\end{align*}
consistent with Crofton \& Exhumatus.

\subsection{\label{PrimalPeri}Primal Perimeter}

As in section [\ref{PrimalCoord}], assume that the triangle $ABC$ satisfies
$\alpha=\theta$, $b=\rho$, $c=\kappa$. These three parameters are related to
perimeter $\tau$ as follows:%

\[
\cos(\theta)=\frac{\sin(\tau-\kappa)}{\sin(\kappa)}+\frac{\cos(\tau
-\kappa)-\cos(\kappa)}{\sin(\kappa)}\cot(\rho).
\]
A\ proof appears in section [\ref{TrigIdentPrimal2}]. For fixed $\tau$ and
$\kappa$, define
\[
f(\rho)=\left\{
\begin{array}
[c]{lll}%
\pi &  & \text{if }0\leq\rho<\tau/2-\kappa,\\
\arccos\left[  \frac{\sin(\tau-\kappa)}{\sin(\kappa)}+\frac{\cos(\tau
-\kappa)-\cos(\kappa)}{\sin(\kappa)}\cot(\rho)\right]  &  & \text{if }%
\tau/2-\kappa\leq\rho\leq\tau/2,\\
0 &  & \text{if }\tau/2<\rho\leq\pi
\end{array}
\right.
\]
assuming $\kappa\leq\tau/2$; otherwise $f(\rho)=0$. Then the conditional
probability is
\begin{align*}
\operatorname*{P}\left\{  \text{perimeter}\leq\tau\left|  c=\kappa\right.
\right\}   &  =\frac1{2\pi}%
{\displaystyle\int\limits_{0}^{\pi}}
{\displaystyle\int\limits_{0}^{f(\rho)}}
\sin(\rho)d\theta\,d\rho\\
&  =\frac1{2\pi}\left(
{\displaystyle\int\limits_{0}^{\tau/2-\kappa}}
{\displaystyle\int\limits_{0}^{\pi}}
\sin(\rho)d\theta\,d\rho+%
{\displaystyle\int\limits_{\tau/2-\kappa}^{\tau/2}}
{\displaystyle\int\limits_{0}^{f(\rho)}}
\sin(\rho)d\theta\,d\rho+%
{\displaystyle\int\limits_{\tau/2}^{\pi}}
{\displaystyle\int\limits_{0}^{0}}
\sin(\rho)d\theta\,d\rho\right) \\
&  =\frac1{2\pi}\left(  \pi\left[  1-\cos\left(  \tfrac\tau2-\kappa\right)
\right]  +%
{\displaystyle\int\limits_{\tau/2-\kappa}^{\tau/2}}
f(\rho)\sin(\rho)d\rho\right)  .
\end{align*}
This result can be experimentally verified by generating many primal triangles
$ABC$ with $c=\kappa$, and then plotting all pairs $(\rho,\theta)$
corresponding to triangles with perimeter $\leq\tau$. The scatterplot fills
the region $[0,\pi]\times[0,\pi]$ except for the portion lying above the curve
$\theta=f(\rho)$.

An exact evaluation of the unconditional probability $\operatorname*{P}%
\left\{  \text{perimeter}\leq\tau\right\}  $, for arbitrary $\tau$, remains
open [\ref{Unsuccess}].

Here is a quick method for evaluating not probability, but instead density, at
$\tau=\pi$. We start with the conditional density
\begin{align*}
&  \ \ \ \ \frac d{d\tau}\frac1{2\pi}\left(  \pi\left[  1-\cos\left(
\tfrac\tau2-\kappa\right)  \right]  +%
{\displaystyle\int\limits_{\tau/2-\kappa}^{\tau/2}}
f(\rho)\sin(\rho)d\rho\right) \\
\  &  =\frac1{2\pi}\left(  \tfrac\pi2\sin\left(  \tfrac\tau2-\kappa\right)
+\tfrac12f\left(  \tfrac\tau2\right)  \sin\left(  \tfrac\tau2\right)
-\tfrac12f\left(  \tfrac\tau2-\kappa\right)  \sin\left(  \tfrac\tau
2-\kappa\right)  +%
{\displaystyle\int\limits_{\tau/2-\kappa}^{\tau/2}}
\frac d{d\tau}f(\rho)\sin(\rho)d\rho\right) \\
\  &  =\frac1{2\pi}%
{\displaystyle\int\limits_{\tau/2-\kappa}^{\tau/2}}
g(\rho)\sin(\rho)d\rho
\end{align*}
where
\[
g(\rho)=\frac{df}{d\tau}=\frac{\sin(\tau-\kappa-\rho)\sin(\rho)}{\sqrt
{\sin(\kappa)^{2}\sin(\rho)^{2}-\left[  \cos(\kappa)\cos(\rho)-\cos
(\tau-\kappa-\rho)\right]  ^{2}}}.
\]
It follows that the unconditional density is
\[
\frac1{4\pi}%
{\displaystyle\int\limits_{0}^{\tau/2}}
{\displaystyle\int\limits_{\tau/2-\kappa}^{\tau/2}}
\frac{\sin(\tau-\kappa-\rho)\sin(\rho)}{\sqrt{\sin(\kappa)^{2}\sin(\rho
)^{2}-\left[  \cos(\kappa)\cos(\rho)-\cos(\tau-\kappa-\rho)\right]  ^{2}}}%
\sin(\kappa)d\rho\,d\kappa
\]
because the density for $\kappa$ is $\sin(\kappa)/2.$ In the special case when
$\tau=\pi$, this becomes
\begin{align*}
&  \ \frac1{4\pi}%
{\displaystyle\int\limits_{0}^{\pi/2}}
{\displaystyle\int\limits_{\pi/2-\kappa}^{\pi/2}}
\frac{\sin(\kappa+\rho)\sin(\kappa)\sin(\rho)}{\sqrt{\sin(\kappa)^{2}\sin
(\rho)^{2}-\left[  \cos(\kappa)\cos(\rho)+\cos(\kappa+\rho)\right]  ^{2}}%
}d\rho\,d\kappa\\
\  &  =\frac1{4\pi}%
{\displaystyle\int\limits_{0}^{\pi/2}}
{\displaystyle\int\limits_{\pi/2-\kappa}^{\pi/2}}
\frac{\sin(\kappa+\rho)\sin(\kappa)\sin(\rho)}{\sqrt{-4\cos(\kappa)\cos
(\rho)\cos(\kappa+\rho)}}d\rho\,d\kappa\\
\  &  =\frac1{8\pi}%
{\displaystyle\int\limits_{0}^{\pi/2}}
{\displaystyle\int\limits_{\pi/2-\kappa}^{\pi/2}}
\frac{\sin(\kappa+\rho)}{\sqrt{-\cos(\kappa+\rho)}}\frac{\sin(\kappa)}%
{\sqrt{\cos(\kappa)}}\frac{\sin(\rho)}{\sqrt{\cos(\rho)}}d\rho\,d\kappa\\
\  &  =\frac1{4\pi}%
{\displaystyle\int\limits_{\pi/4}^{\pi/2}}
{\displaystyle\int\limits_{\pi/2-\kappa}^{\kappa}}
\frac{\sin(\kappa+\rho)}{\sqrt{-\cos(\kappa+\rho)}}\frac{\sin(\kappa)\sin
(\rho)}{\sqrt{\cos(\kappa)\cos(\rho)}}d\rho\,d\kappa\\
\  &  =\frac1{4\sqrt{2}\pi}%
{\displaystyle\int\limits_{\pi/4}^{\pi/2}}
{\displaystyle\int\limits_{\pi/2-\kappa}^{\kappa}}
\frac{\sin(\kappa+\rho)}{\sqrt{-\cos(\kappa+\rho)}}\frac{\cos(\kappa
-\rho)-\cos(\kappa+\rho)}{\sqrt{\cos(\kappa-\rho)+\cos(\kappa+\rho)}}%
d\rho\,d\kappa.
\end{align*}
Let $u=\kappa+\rho$, $v=\kappa-\rho$. Then $\left|  \partial(u,v)/\partial
(\kappa,\rho)\right|  =2$ and the integral is transformed to
\[
\frac1{8\sqrt{2}\pi}%
{\displaystyle\int\limits_{0}^{\pi/2}}
{\displaystyle\int\limits_{\pi/2}^{\pi-v}}
\frac{\sin(u)}{\sqrt{-\cos(u)}}\frac{\cos(v)-\cos(u)}{\sqrt{\cos(v)+\cos(u)}%
}du\,dv=\frac{3\sqrt{2}}{32}
\]
as promised.

\subsection{\label{DualArea}Dual Area}

As in section [\ref{DualCoord}], assume that the triangle $ABC$ satisfies
$\alpha=\kappa$, $\beta=\theta$, $c=\rho$. These three parameters are related
to area $\sigma$ as follows:%

\[
-\cos(\rho)=\frac{\sin(\sigma-\kappa)}{\sin(\kappa)}+\frac{\cos(\sigma
-\kappa)-\cos(\kappa)}{\sin(\kappa)}\cot(\theta).
\]
A\ proof appears in section [\ref{TrigIdentDual2}]. For fixed $\sigma$ and
$\kappa$, define
\[
f(\theta)=\left\{
\begin{array}
[c]{lll}%
\pi &  & \text{if }0\leq\theta<\sigma/2,\\
\pi-\arccos\left[  \frac{\sin(\sigma-\kappa)}{\sin(\kappa)}+\frac{\cos
(\sigma-\kappa)-\cos(\kappa)}{\sin(\kappa)}\cot(\theta)\right]  &  & \text{if
}\sigma/2\leq\theta\leq\pi-(\kappa-\sigma/2),\\
0 &  & \text{if }\pi-(\kappa-\sigma/2)<\theta\leq\pi
\end{array}
\right.
\]
assuming $\kappa\geq\sigma/2$; otherwise $f(\theta)=\pi$. Then the conditional
probability is
\begin{align*}
\operatorname*{P}\left\{  \text{area}\leq\sigma\left|  \alpha=\kappa\right.
\right\}   &  =\frac1{2\pi}%
{\displaystyle\int\limits_{0}^{\pi}}
{\displaystyle\int\limits_{0}^{f(\theta)}}
\sin(\theta)d\rho\,d\theta\\
\  &  =\frac1{2\pi}\left(
{\displaystyle\int\limits_{0}^{\sigma/2}}
{\displaystyle\int\limits_{0}^{\pi}}
\sin(\theta)d\rho\,d\theta+%
{\displaystyle\int\limits_{\sigma/2}^{\pi-(\kappa-\sigma/2)}}
{\displaystyle\int\limits_{0}^{f(\theta)}}
\sin(\theta)d\rho\,d\theta+%
{\displaystyle\int\limits_{\pi-(\kappa-\sigma/2)}^{\pi}}
{\displaystyle\int\limits_{0}^{0}}
\sin(\theta)d\rho\,d\theta\right) \\
\  &  =\frac1{2\pi}\left(  \pi\left[  1-\cos\left(  \tfrac\sigma2\right)
\right]  +%
{\displaystyle\int\limits_{\sigma/2}^{\pi-(\kappa-\sigma/2)}}
f(\theta)\sin(\theta)d\theta\right)  .
\end{align*}
This result can be experimentally verified by generating many dual triangles
$ABC$ with $\alpha=\kappa$, and then plotting all pairs $(\theta,\rho)$
corresponding to triangles with area $\leq\sigma$. The scatterplot fills the
region $[0,\pi]\times[0,\pi]$ except for the portion lying above the curve
$\rho=f(\theta)$.

Since (dual area) = ($2\pi-$ primal perimeter), it is not surprising that
conditional probabilities are so similar.

For completeness' sake, let us compute the conditional density
\begin{align*}
&  \ \ \ \ \ \frac d{d\sigma}\frac1{2\pi}\left(  \pi\left[  1-\cos\left(
\tfrac\sigma2\right)  \right]  +%
{\displaystyle\int\limits_{\sigma/2}^{\pi-(\kappa-\sigma/2)}}
f(\theta)\sin(\theta)d\theta\right) \\
\  &  =\frac1{2\pi}\left(  \tfrac\pi2\sin\left(  \tfrac\sigma2\right)
+\tfrac12f\left(  \pi-\left(  \kappa-\tfrac\sigma2\right)  \right)
\sin\left(  \pi-\left(  \kappa-\tfrac\sigma2\right)  \right)  -\tfrac
12f\left(  \tfrac\sigma2\right)  \sin\left(  \tfrac\sigma2\right)  \right. \\
&  \ \left.  +%
{\displaystyle\int\limits_{\sigma/2}^{\pi-(\kappa-\sigma/2)}}
\frac d{d\sigma}f(\theta)\sin(\theta)d\theta\right) \\
\  &  =\frac1{2\pi}%
{\displaystyle\int\limits_{\sigma/2}^{\pi-(\kappa-\sigma/2)}}
g(\theta)\sin(\theta)d\theta
\end{align*}
where
\[
g(\theta)=\frac{df}{d\sigma}=-\frac{\sin(\sigma-\kappa-\theta)\sin(\theta
)}{\sqrt{\sin(\kappa)^{2}\sin(\theta)^{2}-\left[  \cos(\kappa)\cos
(\theta)-\cos(\sigma-\kappa-\theta)\right]  ^{2}}}.
\]
It follows that the unconditional density is
\[
-\frac1{4\pi}%
{\displaystyle\int\limits_{\sigma/2}^{\pi}}
{\displaystyle\int\limits_{\sigma/2}^{\pi-(\kappa-\sigma/2)}}
\frac{\sin(\sigma-\kappa-\theta)\sin(\theta)}{\sqrt{\sin(\kappa)^{2}%
\sin(\theta)^{2}-\left[  \cos(\kappa)\cos(\theta)-\cos(\sigma-\kappa
-\theta)\right]  ^{2}}}\sin(\kappa)d\theta\,d\kappa
\]
because the density for $\kappa$ is $\sin(\kappa)/2.$ In the special case when
$\sigma=\pi$, this becomes
\begin{align*}
&  \ -\frac1{4\pi}%
{\displaystyle\int\limits_{\pi/2}^{\pi}}
{\displaystyle\int\limits_{\pi/2}^{3\pi/2-\kappa}}
\frac{\sin(\kappa+\theta)\sin(\kappa)\sin(\theta)}{\sqrt{\sin(\kappa)^{2}%
\sin(\theta)^{2}-\left[  \cos(\kappa)\cos(\theta)+\cos(\kappa+\theta)\right]
^{2}}}d\theta\,d\kappa\\
\  &  =-\frac1{4\pi}%
{\displaystyle\int\limits_{\pi/2}^{\pi}}
{\displaystyle\int\limits_{\pi/2}^{3\pi/2-\kappa}}
\frac{\sin(\kappa+\theta)\sin(\kappa)\sin(\theta)}{\sqrt{-4\cos(\kappa
)\cos(\theta)\cos(\kappa+\theta)}}d\theta\,d\kappa\\
\  &  =-\frac1{8\pi}%
{\displaystyle\int\limits_{\pi/2}^{\pi}}
{\displaystyle\int\limits_{\pi/2}^{3\pi/2-\kappa}}
\frac{\sin(\kappa+\theta)}{\sqrt{-\cos(\kappa+\theta)}}\frac{\sin(\kappa
)}{\sqrt{-\cos(\kappa)}}\frac{\sin(\theta)}{\sqrt{-\cos(\theta)}}%
d\theta\,d\kappa\\
\  &  =-\frac1{4\pi}%
{\displaystyle\int\limits_{\pi/2}^{3\pi/4}}
{\displaystyle\int\limits_{\kappa}^{3\pi/2-\kappa}}
\frac{\sin(\kappa+\theta)}{\sqrt{-\cos(\kappa+\theta)}}\frac{\sin(\kappa
)\sin(\theta)}{\sqrt{\cos(\kappa)\cos(\theta)}}d\theta\,d\kappa\\
\  &  =-\frac1{4\sqrt{2}\pi}%
{\displaystyle\int\limits_{\pi/2}^{3\pi/4}}
{\displaystyle\int\limits_{\kappa}^{3\pi/2-\kappa}}
\frac{\sin(\kappa+\theta)}{\sqrt{-\cos(\kappa+\theta)}}\frac{\cos
(\kappa-\theta)-\cos(\kappa+\theta)}{\sqrt{\cos(\kappa-\theta)+\cos
(\kappa+\theta)}}d\theta\,d\kappa\\
\  &  =-\frac1{8\sqrt{2}\pi}%
{\displaystyle\int\limits_{0}^{\pi/2}}
{\displaystyle\int\limits_{\pi+v}^{3\pi/2}}
\frac{\sin(u)}{\sqrt{-\cos(u)}}\frac{\cos(v)-\cos(u)}{\sqrt{\cos(v)+\cos(u)}%
}du\,dv=\frac{3\sqrt{2}}{32}%
\end{align*}
as promised.

\section{Two Evaluations}

\subsection{\label{Success}Successful Evaluation for Primal Area}

Starting from the half-angle formula for tangent
\[
\tan\left(  \frac\rho2\right)  ^{2}=\frac{1-\cos(\rho)}{1+\cos(\rho)}
\]
we deduce that
\[
\cot\left(  \frac\rho2\right)  ^{2}+1=\frac{1+\cos(\rho)}{1-\cos(\rho
)}+1=\frac2{1-\cos(\rho)}
\]
hence
\[
\frac{1-\cos(\rho)}2=\frac1{\cot\left(  \frac\rho2\right)  ^{2}+1}%
=\frac1{\Omega^{2}\sin\left(  \theta-\frac\sigma2\right)  ^{2}+1}
\]
by [\ref{PrimalArea}], where
\[
\Omega=\tan\left(  \frac\kappa2\right)  \csc\left(  \frac\sigma2\right)  .
\]
It follows that
\begin{align*}%
{\displaystyle\int\limits_{\sigma/2}^{\pi}}
\frac{1-\cos(f(\theta))}2d\theta &  =%
{\displaystyle\int\limits_{\sigma/2}^{\pi}}
\frac1{\Omega^{2}\sin\left(  \theta-\frac\sigma2\right)  ^{2}+1}d\theta\\
\  &  =\left\{
\begin{array}
[c]{lll}%
\dfrac{\pi-\arctan\left(  \sqrt{\Omega^{2}+1}\tan(\sigma/2)\right)  }%
{\sqrt{\Omega^{2}+1}} &  & \text{if }0\leq\sigma<\pi,\\
-\dfrac{\arctan\left(  \sqrt{\Omega^{2}+1}\tan(\sigma/2)\right)  }%
{\sqrt{\Omega^{2}+1}} &  & \text{if }\pi\leq\sigma\leq2\pi
\end{array}
\right.
\end{align*}
and thus
\[
\operatorname*{P}\left\{  \text{area}\leq\sigma\right\}  =\left\{
\begin{array}
[c]{lll}%
\dfrac1{2\pi}%
{\displaystyle\int\limits_{0}^{\pi}}
\left[  \dfrac\sigma2+\dfrac{\pi-\arctan\left(  \sqrt{\Omega^{2}+1}\tan
(\sigma/2)\right)  }{\sqrt{\Omega^{2}+1}}\right]  \sin(\kappa)d\kappa &  &
\text{if }0\leq\sigma<\pi,\\
\dfrac1{2\pi}%
{\displaystyle\int\limits_{0}^{\pi}}
\left[  \dfrac\sigma2-\dfrac{\arctan\left(  \sqrt{\Omega^{2}+1}\tan
(\sigma/2)\right)  }{\sqrt{\Omega^{2}+1}}\right]  \sin(\kappa)d\kappa &  &
\text{if }\pi\leq\sigma\leq2\pi.
\end{array}
\right.
\]
The area density is therefore
\[
\left\{
\begin{array}
[c]{lll}%
\dfrac1{2\pi}\left[  1+\dfrac d{d\sigma}%
{\displaystyle\int\limits_{0}^{\pi}}
\dfrac{\pi-\arctan\left(  \sqrt{\Omega^{2}+1}\tan(\sigma/2)\right)  }%
{\sqrt{\Omega^{2}+1}}\sin(\kappa)d\kappa\right]  &  & \text{if }0\leq
\sigma<\pi,\\
\dfrac1{2\pi}\left[  1-\dfrac d{d\sigma}%
{\displaystyle\int\limits_{0}^{\pi}}
\dfrac{\arctan\left(  \sqrt{\Omega^{2}+1}\tan(\sigma/2)\right)  }{\sqrt
{\Omega^{2}+1}}\sin(\kappa)d\kappa\right]  &  & \text{if }\pi\leq\sigma
\leq2\pi
\end{array}
\right.
\]
which, as outlined in [\ref{CroftonExhum}], gives rise to the
Crofton/Exhumatus expression.

\subsection{\label{Unsuccess}Unsuccessful Evaluation for Primal Perimeter}

There does not seem to be an analogous approach for computing
\[%
{\displaystyle\int\limits_{\tau/2-\kappa}^{\tau/2}}
f(\rho)\sin(\rho)d\rho
\]
from [\ref{PrimalPeri}] in closed-form.\ \ We\ suspect that elliptic integrals
will arise, but have not yet found a method for demonstrating this. See
[\ref{JonesBenyTink}] for more details.

\section{Two More Approaches}

\subsection{\label{MedianArea}Median Area}

Assume that the triangle $ABC$ satisfies $c=\kappa$. A \textbf{median} in
$ABC$ is the great circle drawn from vertex $C$ to the midpoint $P$ of side $c
$. Let $\rho$ denote the spherical distance between $P$ and $C$, and $\theta$
denote the angle between $PB$ and $PC$. These three parameters are related to
primal area $\sigma$ as follows:
\[
\tan\left(  \frac\sigma2\right)  =\frac{\sin(\kappa/2)\sin(\rho)\sin(\theta
)}{\cos(\kappa/2)+\cos(\rho)}.
\]
A\ proof appears in section [\ref{TrigMedian}]. For fixed $\sigma$ and
$\kappa$, define
\[
f(\theta)=\arccos\left[  \frac{\cos\left(  \tfrac\kappa2\right)  \tan\left(
\tfrac\sigma2\right)  ^{2}\csc(\theta)^{2}\mp\sin\left(  \tfrac\kappa2\right)
^{2}\sqrt{1+\tan\left(  \tfrac\sigma2\right)  ^{2}\csc(\theta)^{2}}}%
{\cos\left(  \tfrac\kappa2\right)  ^{2}-\tan\left(  \tfrac\sigma2\right)
^{2}\csc(\theta)^{2}-1}\right]
\]
where $-$ is chosen if $\sigma\leq\pi$ and $+$ is chosen if $\sigma>\pi$; then
the conditional probability is
\begin{align*}
\operatorname*{P}\left\{  \text{area}\leq\sigma\left|  c=\kappa\right.
\right\}   &  =\frac1{2\pi}%
{\displaystyle\int\limits_{0}^{\pi}}
{\displaystyle\int\limits_{0}^{f(\theta)}}
\sin(\rho)d\rho\,d\theta\\
\  &  =\frac1{2\pi}%
{\displaystyle\int\limits_{0}^{\pi}}
(1-\cos(f(\theta)))d\theta\\
\  &  =\dfrac1{2\pi}%
{\displaystyle\int\limits_{0}^{\pi}}
\left(  1-\frac{\cos\left(  \tfrac\kappa2\right)  \tan\left(  \tfrac
\sigma2\right)  ^{2}\csc(\theta)^{2}\mp\sin\left(  \tfrac\kappa2\right)
^{2}\sqrt{1+\tan\left(  \tfrac\sigma2\right)  ^{2}\csc(\theta)^{2}}}%
{\cos\left(  \tfrac\kappa2\right)  ^{2}-\tan\left(  \tfrac\sigma2\right)
^{2}\csc(\theta)^{2}-1}\right)  d\theta.
\end{align*}
This result can be experimentally verified by generating many primal triangles
$ABC$ with $c=\kappa$, and then plotting all pairs $(\theta,\rho)$
corresponding to triangles with area $\leq\sigma$. The scatterplot fills the
region $[0,\pi]\times[0,\pi]$ except for the portion lying above the curve
$\rho=f(\theta)$. This approach is believed to be the same as Crofton \&
Exhumatus (details in \cite{CE} are rather thin). Not seeing any advantage
over our approach in [\ref{PrimalArea}], we stop here.

\subsection{\label{BisectorPerimeter}Bisector Perimeter}

Assume that the triangle $ABC$ satisfies $\alpha=\kappa$. An \textbf{angle
bisector} in $ABC$ is the great circle drawn from vertex $A$ that splits angle
$\alpha$ in half. Define $Q$ to be the intersection point between this circle
and side $BC$. Let $\rho$ denote the spherical distance between $Q$ and $A$,
and $\theta$ denote the angle between $QC$ and $QA$. These three parameters
are related to dual perimeter $\tau$ as follows:
\[
\tan\left(  \frac\tau2\right)  =-\frac{\cos(\kappa/2)\sin(\rho)\sin(\theta
)}{\sin(\kappa/2)+\cos(\rho)\sin(\theta)}.
\]
A\ proof appears in section [\ref{TrigBisector}]. For fixed $\tau$ and
$\kappa$, define
\[
f_{\text{lower}}(\rho)=\arcsin\left[  -\frac{\tan(\tau/2)\sin(\kappa/2)}%
{\tan(\tau/2)\cos(\rho)+\cos(\kappa/2)\sin(\rho)}\right]  ,
\]
\[
f_{\text{upper}}(\rho)=\pi-\arcsin\left[  -\frac{\tan(\tau/2)\sin(\kappa
/2)}{\tan(\tau/2)\cos(\rho)+\cos(\kappa/2)\sin(\rho)}\right]
\]
assuming
\[
\arccos\left(  -\frac{\cos(\tau/2)+\sin(\kappa/2)}{1+\cos(\tau/2)\sin
(\kappa/2)}\right)  =\rho_{\text{thres}}\leq\rho\leq\pi.
\]
The region of all pairs $(\rho,\theta)$ corresponding to triangles with
perimeter $\leq\tau$ is more complicated than earlier examples. The
scatterplot fills the region $[\rho_{\text{thres}},\pi]\times[0,\pi]$ except
for portions lying either above the curve $\theta=f_{\text{upper}}(\rho)$ or
below the curve $\theta=f_{\text{lower}}(\rho)$. The conditional probability
is
\begin{align*}
\operatorname*{P}\left\{  \text{perimeter}\leq\tau\left|  \alpha
=\kappa\right.  \right\}   &  =\frac1{2\pi}%
{\displaystyle\int\limits_{\rho_{\text{thres}}}^{\pi}}
{\displaystyle\int\limits_{f_{\text{lower}}(\rho)}^{f_{\text{upper}}(\rho)}}
\sin(\theta)d\theta\,d\rho\\
\  &  =\frac1{2\pi}%
{\displaystyle\int\limits_{\rho_{\text{thres}}}^{\pi}}
(\cos(f_{\text{lower}}(\rho))-\cos(f_{\text{upper}}(\rho)))d\rho\\
\  &  =\frac1\pi%
{\displaystyle\int\limits_{\rho_{\text{thres}}}^{\pi}}
\sqrt{1-\left(  \frac{\tan(\tau/2)\sin(\kappa/2)}{\tan(\tau/2)\cos(\rho
)+\cos(\kappa/2)\sin(\rho)}\right)  ^{2}}d\rho.
\end{align*}
Due to the unanticipated complexity, we stop here.

\section{Two More Coordinate Systems}

The primal coordinate system [\ref{PrimalCoord}] allows us to specify a
triangle, given a fixed side $\kappa$, with an additional side $\rho$ and an
angle $\theta$. Can we do as well with two angles instead? The dual coordinate
system [\ref{DualCoord}] allows us to likewise specify a triangle, given a
fixed angle $\kappa$. Can we do as well with two sides instead?

\subsection{\label{AngleCoord}Angle Coordinates}

Without loss of generality, let $A=(1,0,0)$ and $B=(\cos(\kappa),\sin
(\kappa),0)$ in $xyz$ coordinates. We wish to locate the unique point $C$ in
the hemisphere so that the triangle $ABC$ satisfies $\alpha=\varphi$,
$\beta=\psi$, $c=\kappa$. See Figure 3.
\begin{figure}[ptb]%
\centering
\includegraphics[
height=3.5146in,
width=3.4956in
]%
{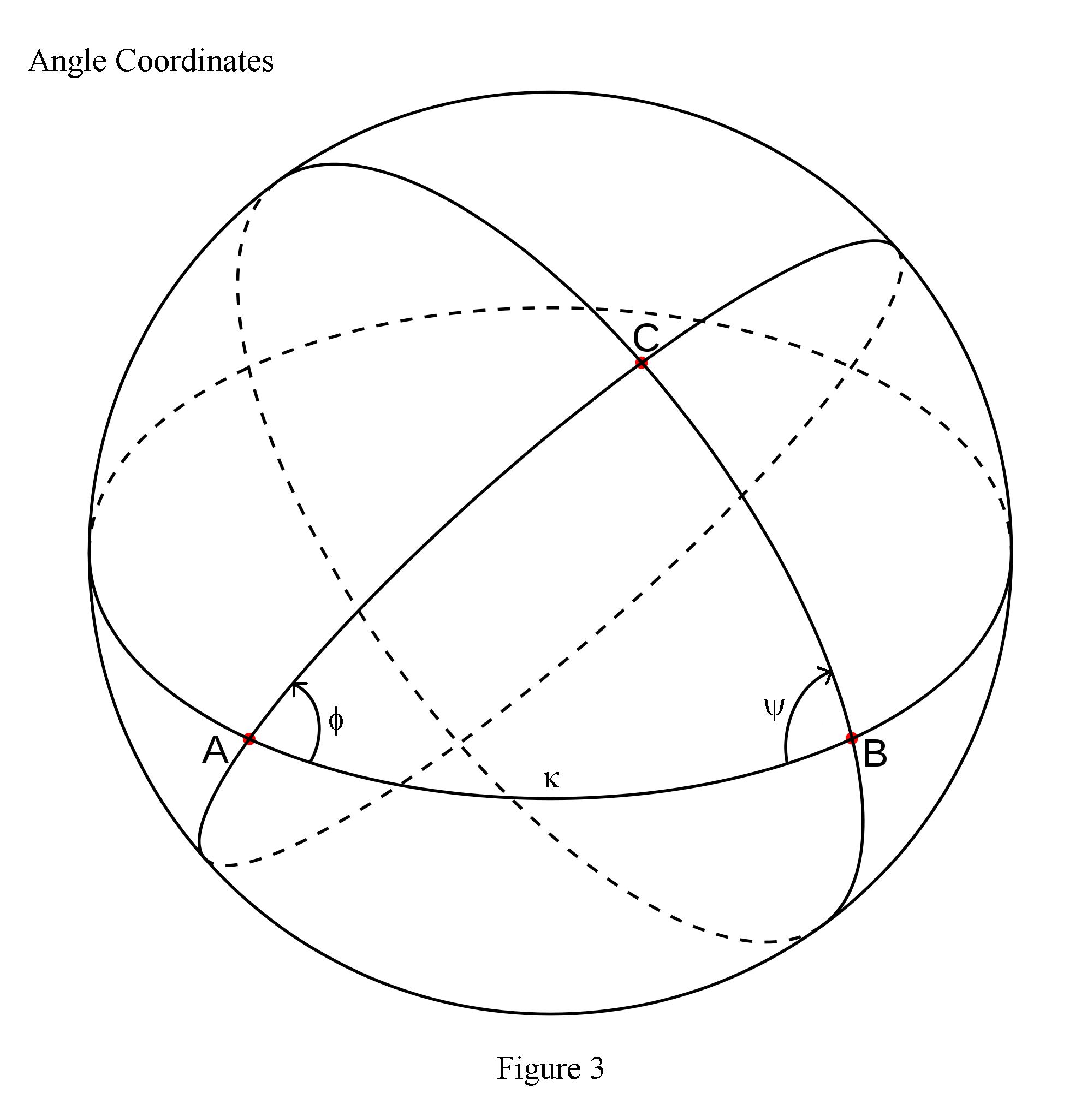}%
\end{figure}
The parameters $\varphi$, $\psi$ are regarded as varying while the parameter
$\kappa$ is fixed.

Let us examine the great circle containing $A$, $C$. It must also contain the
point $(0,\cos(\varphi),\sin(\varphi))$, since this is the image of $(0,1,0)$
after rotation through angle $\varphi$. Hence a normal vector is
$V=(1,0,0)\times(0,\cos(\varphi),\sin(\varphi))=(0,-\sin(\varphi),\cos
(\varphi))$.

Let us examine the great circle containing $B$, $C$. Think of rotating the
equatorial disk in space so that the vector $(\cos(\kappa),\sin(\kappa),0)$
remains fixed and the vector $(\sin(\kappa),-\cos(\kappa),0)$ moves toward
$(0,0,1)$ through the angle $\psi$. The rotation matrix performing this motion
is \cite{Wi}
\[
S=\left(
\begin{array}
[c]{ccc}%
\cos(\kappa)^{2}+(1-\cos(\kappa)^{2})\cos(\psi) & \cos(\kappa)\sin
(\kappa)(1-\cos(\psi)) & -\sin(\kappa)\sin(\psi)\\
\cos(\kappa)\sin(\kappa)(1-\cos(\psi)) & \sin(\kappa)^{2}+(1-\sin(\kappa
)^{2})\cos(\psi) & \cos(\kappa)\sin(\psi)\\
\sin(\kappa)\sin(\psi) & -\cos(\kappa)\sin(\psi) & \cos(\psi)
\end{array}
\right)
\]
and
\[
S\left(
\begin{array}
[c]{c}%
\sin(\kappa)\\
-\cos(\kappa)\\
0
\end{array}
\right)  =\left(
\begin{array}
[c]{c}%
\sin(\kappa)\cos(\psi)\\
-\cos(\kappa)\cos(\psi)\\
\sin(\psi)
\end{array}
\right)  .
\]
Hence the great circle must contain the point $(\sin(\kappa)\cos(\psi
),-\cos(\kappa)\cos(\psi),\sin(\psi))$ and a normal vector is
\[
W=\left(
\begin{array}
[c]{c}%
\cos(\kappa)\\
\sin(\kappa)\\
0
\end{array}
\right)  \times\left(
\begin{array}
[c]{c}%
\sin(\kappa)\cos(\psi)\\
-\cos(\kappa)\cos(\psi)\\
\sin(\psi)
\end{array}
\right)  =\left(
\begin{array}
[c]{c}%
\sin(\kappa)\sin(\psi)\\
-\cos(\kappa)\sin(\psi)\\
-\cos(\psi)
\end{array}
\right)  .
\]
The point $C$ is orthogonal to the two normal vectors and at unit distance
from the origin, equivalently, $C=(V\times W)/\left\|  V\times W\right\|  $.
We have
\[
V\times W=\left(
\begin{array}
[c]{c}%
\sin(\varphi)\cos(\psi)+\cos(\kappa)\cos(\varphi)\sin(\psi)\\
\sin(\kappa)\cos(\varphi)\sin(\psi)\\
\sin(\kappa)\sin(\varphi)\sin(\psi)
\end{array}
\right)  ,
\]
\[
\left\|  V\times W\right\|  =\sqrt{1-\left(  \cos(\varphi)\cos(\psi
)-\cos(\kappa)\sin(\varphi)\sin(\psi)\right)  ^{2}}
\]
and thus the Jacobian determinant of $(r,\varphi,\psi)\mapsto(x,y,z)=-rC$
simplifies to
\[
\frac{\sin(\kappa)^{2}\sin(\varphi)\sin(\psi)\left[  \left(  \sin(\varphi
)\cos(\psi)+\cos(\kappa)\cos(\varphi)\sin(\psi)\right)  ^{2}+\sin(\kappa
)^{2}\sin(\psi)^{2}\right]  }{\left[  1-\left(  \cos(\varphi)\cos(\psi
)-\cos(\kappa)\sin(\varphi)\sin(\psi)\right)  ^{2}\right]  ^{5/2}}.
\]

\subsection{\label{SideCoord}Side Coordinates}

Without loss of generality, let $A=(1,0,0)$ and $B=(\cos(\xi),\sin(\xi),0)$ in
$xyz$ coordinates. It seems (at first glance) that we should locate the unique
point $C$ in the upper hemisphere so that the triangle $ABC$ satisfies
$\alpha=\kappa$, $c=\xi$, $b=\eta$. See Figure 4.
\begin{figure}[ptb]%
\centering
\includegraphics[
height=3.589in,
width=3.5682in
]%
{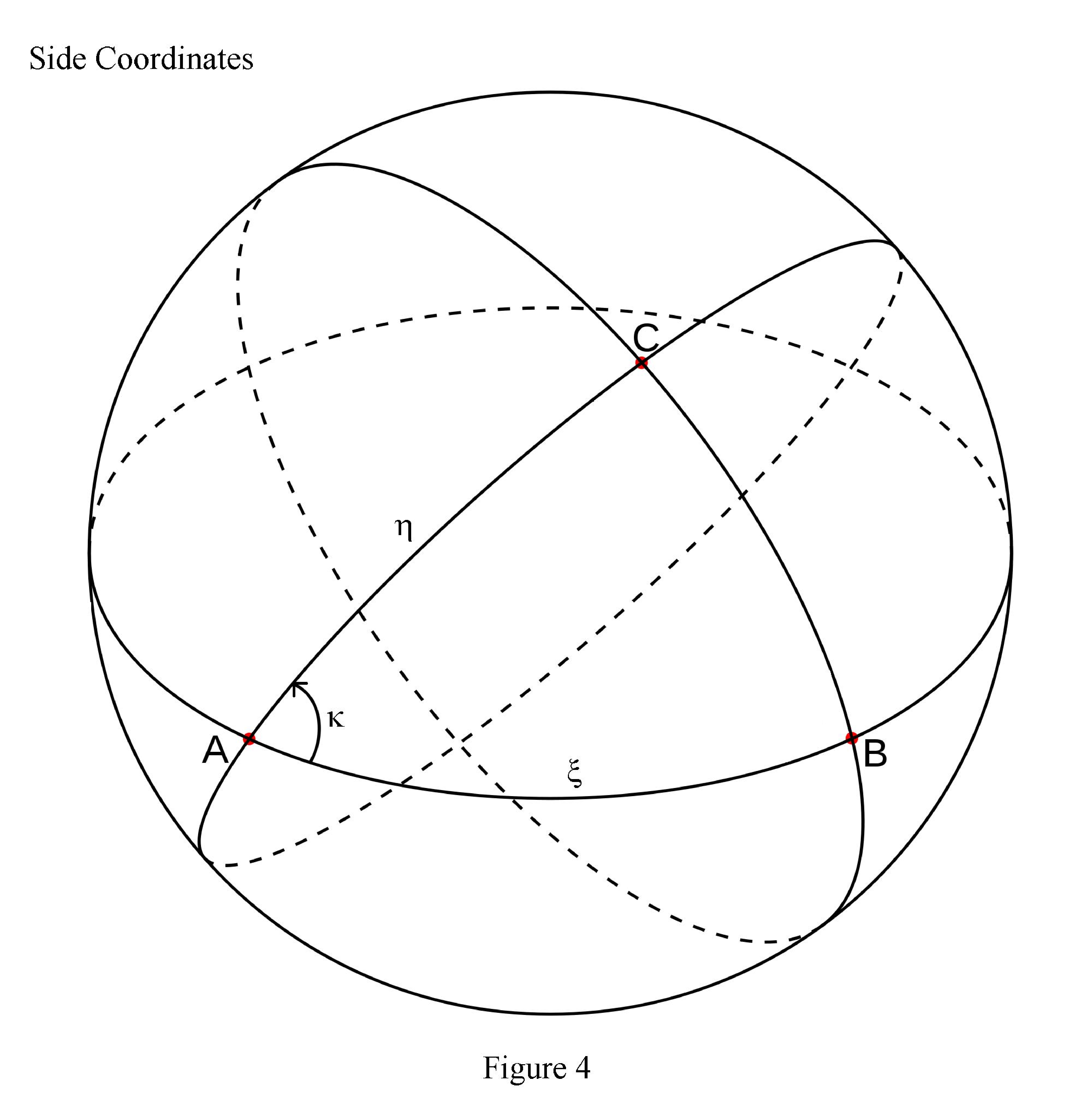}%
\end{figure}
The parameters $\xi$, $\eta$ are regarded as varying while the parameter
$\kappa$ is fixed. Think of rotating the equatorial disk in space so that the
vector $(1,0,0)$ remains fixed and the vector $(0,1,0)$ moves toward $(0,0,1)$
through the angle $\kappa$. The rotation matrix performing this motion is
\cite{Wi}
\[
R=\left(
\begin{array}
[c]{ccc}%
1 & 0 & 0\\
0 & \cos(\kappa) & -\sin(\kappa)\\
0 & \sin(\kappa) & \cos(\kappa)
\end{array}
\right)
\]
and
\[
R\left(
\begin{array}
[c]{c}%
\cos(\eta)\\
\sin(\eta)\\
0
\end{array}
\right)  =\left(
\begin{array}
[c]{c}%
\cos(\eta)\\
\cos(\kappa)\sin(\eta)\\
\sin(\kappa)\sin(\eta)
\end{array}
\right)  ,
\]
which gives the point $C$.

Now, in fact, this is less than what is required. We must (on second glance)
specify the great circle containing $B$, $C$. This is done via a normal vector
$U=(B\times C)/\left\|  B\times C\right\|  $, where
\begin{align*}
B\times C  &  =\left(
\begin{array}
[c]{c}%
\cos(\xi)\\
\sin(\xi)\\
0
\end{array}
\right)  \times\left(
\begin{array}
[c]{c}%
\cos(\eta)\\
\cos(\kappa)\sin(\eta)\\
\sin(\kappa)\sin(\eta)
\end{array}
\right) \\
\  &  =\left(
\begin{array}
[c]{c}%
\sin(\kappa)\sin(\xi)\sin(\eta)\\
-\sin(\kappa)\cos(\xi)\sin(\eta)\\
-\sin(\xi)\cos(\eta)+\cos(\kappa)\cos(\xi)\sin(\eta)
\end{array}
\right)  ,
\end{align*}
\[
\left\|  B\times C\right\|  =\sqrt{1-\left(  \cos(\xi)\cos(\eta)+\cos
(\kappa)\sin(\xi)\sin(\eta)\right)  ^{2}}.
\]
Thus the Jacobian determinant of $(r,\xi,\eta)\mapsto(x,y,z)=rU$ simplifies
to
\[
\frac{\sin(\kappa)^{2}\sin(\xi)\sin(\eta)\left[  \left(  \sin(\xi)\cos
(\eta)-\cos(\kappa)\cos(\xi)\sin(\eta)\right)  ^{2}+\sin(\kappa)^{2}\sin
(\eta)^{2}\right]  }{\left[  1-\left(  \cos(\xi)\cos(\eta)+\cos(\kappa
)\sin(\xi)\sin(\eta)\right)  ^{2}\right]  ^{5/2}}.
\]

\subsection{Possible Applications}

Combining an identity in \cite{Er} with the Law of Cosines for Angles, we
obtain
\[
\tan\left(  \frac{\tau}{2}\right)  =\frac{\sin(\varphi)\sin(\psi)\sin(\kappa
)}{\cos(\varphi)+\cos(\psi)-\cos(\varphi)\cos(\psi)+\cos(\kappa)\sin
(\varphi)\sin(\psi)-1}.
\]
Let us solve for $\psi$ as follows:
\[
\left[  \cos(\varphi)+\cos(\psi)-\cos(\varphi)\cos(\psi)-1\right]  \tan\left(
\frac{\tau}{2}\right)  =\sin(\varphi)\sin(\psi)\left[  \sin(\kappa
)-\cos(\kappa)\tan\left(  \frac{\tau}{2}\right)  \right]
\]
hence
\[
-\left[  1-\cos(\varphi)\right]  \left[  1-\cos(\psi)\right]  \tan\left(
\frac{\tau}{2}\right)  =\sin(\varphi)\sin(\psi)\left[  \sin(\kappa
)-\cos(\kappa)\tan\left(  \frac{\tau}{2}\right)  \right]
\]
hence
\[
-\frac{1-\cos(\varphi)}{\sin(\varphi)}\tan\left(  \frac{\tau}{2}\right)
=\frac{\sin(\psi)}{1-\cos(\psi)}\left[  \sin(\kappa)-\cos(\kappa)\tan\left(
\frac{\tau}{2}\right)  \right]
\]
hence
\[
\tan\left(  \frac{\varphi}{2}\right)  \csc\left(  \frac{\tau}{2}%
-\kappa\right)  \sin\left(  \frac{\tau}{2}\right)  =\frac{\sin(\psi)}%
{1-\cos(\psi)}%
\]
hence
\[
\cos(\psi)=\frac{\tan\left(  \frac{\varphi}{2}\right)  ^{2}\csc\left(
\frac{\tau}{2}-\kappa\right)  ^{2}\sin\left(  \frac{\tau}{2}\right)  ^{2}%
-1}{\tan\left(  \frac{\varphi}{2}\right)  ^{2}\csc\left(  \frac{\tau}%
{2}-\kappa\right)  ^{2}\sin\left(  \frac{\tau}{2}\right)  ^{2}+1}%
\]
because $y=\sqrt{1-x^{2}}/(1-x)$ has inverse $x=\left(  y^{2}-1\right)
/\left(  y^{2}+1\right)  $. For fixed $\tau$ and $\kappa$, define
\[
f(\varphi)=\left\{
\begin{array}
[c]{lll}%
\arccos\left[  \frac{\tan\left(  \frac{\varphi}{2}\right)  ^{2}\csc\left(
\frac{\tau}{2}-\kappa\right)  ^{2}\sin\left(  \frac{\tau}{2}\right)  ^{2}%
-1}{\tan\left(  \frac{\varphi}{2}\right)  ^{2}\csc\left(  \frac{\tau}%
{2}-\kappa\right)  ^{2}\sin\left(  \frac{\tau}{2}\right)  ^{2}+1}\right]  &  &
\text{if }0\leq\kappa<\tau/2,\\
0 &  & \text{if }\tau/2\leq\kappa\leq\pi
\end{array}
\right.
\]
then the conditional probability, given $c$, is
\[
\operatorname*{P}\left\{  \text{perimeter}\leq\tau\left\vert c=\kappa\right.
\right\}  =\frac{1}{2\pi}%
{\displaystyle\int\limits_{0}^{\pi}}
{\displaystyle\int\limits_{0}^{f(\varphi)}}
\tfrac{\sin(\kappa)^{2}\sin(\varphi)\sin(\psi)\left[  \left(  \sin
(\varphi)\cos(\psi)+\cos(\kappa)\cos(\varphi)\sin(\psi)\right)  ^{2}%
+\sin(\kappa)^{2}\sin(\psi)^{2}\right]  }{\left[  1-\left(  \cos(\varphi
)\cos(\psi)-\cos(\kappa)\sin(\varphi)\sin(\psi)\right)  ^{2}\right]  ^{5/2}%
}d\psi\,d\varphi.
\]

Similarly, combining an identity in \cite{Er} with the Law of Cosines for
Sides, we obtain
\[
\tan\left(  \frac\sigma2\right)  =\frac{\sin(\xi)\sin(\eta)\sin(\kappa
)}{1+\cos(\xi)+\cos(\eta)+\cos(\xi)\cos(\eta)+\cos(\kappa)\sin(\xi)\sin(\eta
)}.
\]
Let us solve for $\eta$ as follows:
\[
\left[  1+\cos(\xi)+\cos(\eta)+\cos(\xi)\cos(\eta)\right]  \tan\left(
\frac\sigma2\right)  =\sin(\xi)\sin(\eta)\left[  \sin(\kappa)-\cos(\kappa
)\tan\left(  \frac\sigma2\right)  \right]
\]
hence
\[
\frac{1+\cos(\xi)}{\sin(\xi)}\tan\left(  \frac\sigma2\right)  =\frac{\sin
(\eta)}{1+\cos(\eta)}\left[  \sin(\kappa)-\cos(\kappa)\tan\left(  \frac
\sigma2\right)  \right]
\]
hence
\[
\cos(\eta)=\frac{1-\cot\left(  \frac\xi2\right)  ^{2}\csc\left(  \kappa
-\frac\sigma2\right)  ^{2}\sin\left(  \frac\sigma2\right)  ^{2}}{1+\cot\left(
\frac\xi2\right)  ^{2}\csc\left(  \kappa-\frac\sigma2\right)  ^{2}\sin\left(
\frac\sigma2\right)  ^{2}}
\]
because $y=\sqrt{1-x^{2}}/(1+x)$ has inverse $x=\left(  1-y^{2}\right)
/\left(  1+y^{2}\right)  $. For fixed $\sigma$ and $\kappa$, define
\[
f(\xi)=\left\{
\begin{array}
[c]{lll}%
\pi &  & \text{if }0\leq\kappa<\sigma/2,\\
\arccos\left[  \frac{1-\cot\left(  \frac\xi2\right)  ^{2}\csc\left(
\kappa-\frac\sigma2\right)  ^{2}\sin\left(  \frac\sigma2\right)  ^{2}}%
{1+\cot\left(  \frac\xi2\right)  ^{2}\csc\left(  \kappa-\frac\sigma2\right)
^{2}\sin\left(  \frac\sigma2\right)  ^{2}}\right]  &  & \text{if }\sigma
/2\leq\kappa\leq\pi
\end{array}
\right.
\]
then the conditional probability, given $\alpha$, is
\begin{align*}
\operatorname*{P}\left\{  \text{area}\leq\sigma\left|  \alpha=\kappa\right.
\right\}  =\frac1{2\pi}%
{\displaystyle\int\limits_{0}^{\pi}}
{\displaystyle\int\limits_{0}^{f(\xi)}}
\tfrac{\sin(\kappa)^{2}\sin(\xi)\sin(\eta)\left[  \left(  \sin(\xi)\cos
(\eta)-\cos(\kappa)\cos(\xi)\sin(\eta)\right)  ^{2}+\sin(\kappa)^{2}\sin
(\eta)^{2}\right]  }{\left[  1-\left(  \cos(\xi)\cos(\eta)+\cos(\kappa
)\sin(\xi)\sin(\eta)\right)  ^{2}\right]  ^{5/2}}d\eta\,d\xi.
\end{align*}
We have not further pursued this direction of inquiry.

\section{Trigonometric Identities}

The following formulas are required in the main text.

\subsection{\label{TrigIdentPrimal1}Primal Case i}

To prove
\[
\tan\left(  \frac b2\right)  =\cot\left(  \frac c2\right)  \csc\left(
\alpha-\frac\sigma2\right)  \sin\left(  \frac\sigma2\right)
\]
we expand $\cos(\sigma/2)$ and make use of Delambre's analogies \cite{Td}:%

\begin{align*}
\cos\left(  \frac\sigma2\right)   &  =\cos\left(  \frac{\beta+\gamma}%
2-\frac{\pi-\alpha}2\right) \\
\  &  =\cos\left(  \frac{\beta+\gamma}2\right)  \sin\left(  \frac
\alpha2\right)  +\sin\left(  \frac{\beta+\gamma}2\right)  \cos\left(
\frac\alpha2\right) \\
\  &  =\left[  \cos\left(  \frac{b+c}2\right)  \sin\left(  \frac
\alpha2\right)  ^{2}+\cos\left(  \frac{b-c}2\right)  \cos\left(  \frac
\alpha2\right)  ^{2}\right]  \sec\left(  \frac a2\right) \\
\  &  =\left[  \left(  \cos\frac b2\cos\frac c2-\sin\frac b2\sin\frac
c2\right)  \frac{1-\cos\alpha}2+\left(  \cos\frac b2\cos\frac c2+\sin\frac
b2\sin\frac c2\right)  \frac{1+\cos\alpha}2\right]  \sec\frac a2\\
\  &  =\left(  \cos\frac b2\cos\frac c2+\sin\frac b2\sin\frac c2\cos
\alpha\right)  \sec\frac a2.
\end{align*}
Also
\begin{align*}
\sin\left(  \frac\sigma2\right)   &  =\sin\left(  \frac{\beta+\gamma}%
2-\frac{\pi-\alpha}2\right) \\
\  &  =-\cos\left(  \frac{\beta+\gamma}2\right)  \cos\left(  \frac
\alpha2\right)  +\sin\left(  \frac{\beta+\gamma}2\right)  \sin\left(
\frac\alpha2\right) \\
\  &  =\left[  -\cos\left(  \frac{b+c}2\right)  \sin\left(  \frac
\alpha2\right)  \cos\left(  \frac\alpha2\right)  +\cos\left(  \frac
{b-c}2\right)  \cos\left(  \frac\alpha2\right)  \sin\left(  \frac
\alpha2\right)  \right]  \sec\left(  \frac a2\right) \\
\  &  =\left[  -\left(  \cos\frac b2\cos\frac c2-\sin\frac b2\sin\frac
c2\right)  +\left(  \cos\frac b2\cos\frac c2+\sin\frac b2\sin\frac c2\right)
\right]  \cos\frac\alpha2\sin\frac\alpha2\sec\frac a2\\
\  &  =2\sin\frac b2\sin\frac c2\cos\frac\alpha2\sin\frac\alpha2\sec\frac
a2=\sin\frac b2\sin\frac c2\sin\alpha\sec\frac a2.
\end{align*}
Dividing, we obtain
\[
\cot\frac\sigma2=\frac{\cos\frac b2\cos\frac c2+\sin\frac b2\sin\frac
c2\cos\alpha}{\sin\frac b2\sin\frac c2\sin\alpha}
\]
hence
\[
\cos\frac\sigma2\sin\alpha=\left(  \cot\frac b2\cot\frac c2+\cos\alpha\right)
\sin\frac\sigma2
\]
and therefore
\[
\sin\left(  \alpha-\frac\sigma2\right)  =\sin\alpha\cos\frac\sigma2-\cos
\alpha\sin\frac\sigma2=\cot\frac b2\cot\frac c2\sin\frac\sigma2
\]
as was to be shown.

\subsection{\label{TrigIdentPrimal2}Primal Case ii}

To prove
\[
\cos(\alpha)=\frac{\sin(\tau-c)}{\sin(c)}+\frac{\cos(\tau-c)-\cos(c)}{\sin
(c)}\cot(b)
\]
we expand $\cos(a)$:
\begin{align*}
\cos(a)  &  =\cos(\tau-b-c)\\
\  &  =\cos\left(  -b+(\tau-c)\right) \\
\  &  =\cos(b)\cos(\tau-c)+\sin(b)\sin(\tau-c)
\end{align*}
and make use of the Law of Cosines for Sides:
\[
\cos(a)=\cos(b)\cos(c)+\sin(b)\sin(c)\cos(\alpha)
\]
thus
\[
\sin(b)\sin(c)\cos(\alpha)=\sin(\tau-c)\sin(b)+\left[  \cos(\tau
-c)-\cos(c)\right]  \cos(b).
\]
Alternatively, we have
\[
\sin(c)\cos(\alpha)=\sin(\tau-c)+\left(  \cos(\tau-c)-\cos(c)\right)  \cot(b)
\]
hence
\[
-\tan(b)=\frac{\cos(\tau-c)-\cos(c)}{\sin(\tau-c)-\sin(c)\cos(\alpha)}
\]
but solving for $b$ turns out to be more complicated than our strategy of
solving for $\alpha$.

\subsection{\label{TrigIdentDual1}Dual Case i}

To prove
\[
\tan\left(  \frac\beta2\right)  =\cot\left(  \frac\alpha2\right)  \sin\left(
\frac\tau2-c\right)  \csc\left(  \frac\tau2\right)
\]
we expand $\cos(\tau/2)$ and make use of Delambre's analogies \cite{Td}:%

\begin{align*}
\cos\left(  \frac\tau2\right)   &  =\cos\left(  \frac{a+b}2+\frac c2\right) \\
\  &  =\cos\left(  \frac{a+b}2\right)  \cos\left(  \frac c2\right)
-\sin\left(  \frac{a+b}2\right)  \sin\left(  \frac c2\right) \\
\  &  =\left[  \cos\left(  \frac{\alpha+\beta}2\right)  \cos\left(  \frac
c2\right)  ^{2}-\cos\left(  \frac{\alpha-\beta}2\right)  \sin\left(  \frac
c2\right)  ^{2}\right]  \csc\left(  \frac\gamma2\right) \\
\  &  =\left[  \left(  \cos\frac\alpha2\cos\frac\beta2-\sin\frac\alpha
2\sin\frac\beta2\right)  \frac{1+\cos c}2-\left(  \cos\frac\alpha2\cos
\frac\beta2+\sin\frac\alpha2\sin\frac\beta2\right)  \frac{1-\cos c}2\right]
\csc\frac\gamma2\\
\  &  =\left(  -\sin\frac\alpha2\sin\frac\beta2+\cos\frac\alpha2\cos\frac
\beta2\cos c\right)  \csc\frac\gamma2.
\end{align*}
Also
\begin{align*}
\sin\left(  \frac\tau2\right)   &  =\sin\left(  \frac{a+b}2+\frac c2\right) \\
\  &  =\cos\left(  \frac{a+b}2\right)  \sin\left(  \frac c2\right)
+\sin\left(  \frac{a+b}2\right)  \cos\left(  \frac c2\right) \\
\  &  =\left[  \cos\left(  \frac{\alpha+\beta}2\right)  \cos\left(  \frac
c2\right)  \sin\left(  \frac c2\right)  +\cos\left(  \frac{\alpha-\beta
}2\right)  \sin\left(  \frac c2\right)  \cos\left(  \frac c2\right)  \right]
\csc\left(  \frac\gamma2\right) \\
\  &  =\left[  \left(  \cos\frac\alpha2\cos\frac\beta2-\sin\frac\alpha
2\sin\frac\beta2\right)  +\left(  \cos\frac\alpha2\cos\frac\beta2+\sin
\frac\alpha2\sin\frac\beta2\right)  \right]  \cos\frac c2\sin\frac c2\csc
\frac\gamma2\\
\  &  =2\cos\frac\alpha2\cos\frac\beta2\cos\frac c2\sin\frac c2\csc\frac
\gamma2=\cos\frac\alpha2\cos\frac\beta2\sin c\csc\frac\gamma2.
\end{align*}
Dividing, we obtain
\[
\cot\frac\tau2=\frac{-\sin\frac\alpha2\sin\frac\beta2+\cos\frac\alpha
2\cos\frac\beta2\cos c}{\cos\frac\alpha2\cos\frac\beta2\sin c}
\]
hence
\[
\cos\frac\tau2\sin c=\left(  -\tan\frac\alpha2\tan\frac\beta2+\cos c\right)
\sin\frac\tau2
\]
and therefore
\[
\sin\left(  \frac\tau2-c\right)  =\sin\frac\tau2\cos c-\cos\frac\tau2\sin
c=\tan\frac\alpha2\tan\frac\beta2\sin\frac\tau2
\]
as was to be shown.

\subsection{\label{TrigIdentDual2}Dual Case ii}

To prove
\[
-\cos(c)=\frac{\sin(\sigma-\alpha)}{\sin(\alpha)}+\frac{\cos(\sigma
-\alpha)-\cos(\alpha)}{\sin(\alpha)}\cot(\beta)
\]
we expand $-\cos(\gamma)$:
\begin{align*}
-\cos(\gamma)  &  =-\cos(\sigma-\alpha-\beta+\pi)\\
\  &  =\cos\left(  -\beta+(\sigma-\alpha)\right) \\
\  &  =\cos(\beta)\cos(\sigma-\alpha)+\sin(\beta)\sin(\sigma-\alpha)
\end{align*}
and make use of the Law of Cosines for Angles:
\[
-\cos(\gamma)=\cos(\alpha)\cos(\beta)-\sin(\alpha)\sin(\beta)\cos(c)
\]
thus
\[
-\sin(\alpha)\sin(\beta)\cos(c)=\sin(\sigma-\alpha)\sin(\beta)+\left[
\cos(\sigma-\alpha)-\cos(\alpha)\right]  \cos(\beta).
\]
Alternatively, we have
\[
-\sin(\alpha)\cos(c)=\sin(\sigma-\alpha)+\left(  \cos(\sigma-\alpha
)-\cos(\alpha)\right)  \cot(\beta)
\]
hence
\[
-\tan(\beta)=\frac{\cos(\sigma-\alpha)-\cos(\alpha)}{\sin(\sigma-\alpha
)+\sin(\alpha)\cos(c)}
\]
but solving for $\beta$ turns out to be more complicated than our strategy of
solving for $c$.

\subsection{\label{TrigMedian}Median Case}

Let $\rho$, $\theta$ be defined within triangle $ABC$ as in \textbf{[}%
\ref{MedianArea}\textbf{]. }Applying the Law of Cosines for Sides to both
triangles $CPB$ and $CPA$, we have
\begin{equation}
\cos(a)=\cos(\rho)\cos(c/2)+\sin(\rho)\sin(c/2)\cos(\theta), \label{eq61}%
\end{equation}
\begin{equation}
\cos(b)=\cos(\rho)\cos(c/2)-\sin(\rho)\sin(c/2)\cos(\theta) \label{eq62}%
\end{equation}
because $\cos(\pi-\theta)=-\cos(\theta)$; hence
\begin{align*}
1+\cos(a)+\cos(b)+\cos(c)  &  =1+2\cos(\rho)\cos(c/2)+\cos(c)\\
\  &  =2\cos(\rho)\cos(c/2)+2\cos(c/2)^{2}\\
\  &  =2\cos(c/2)\left(  \cos(\rho)+\cos(c/2)\right)  .
\end{align*}
Applying the Law of Sines to both triangles $CPB$ and $ABC$, we have
\[%
\begin{array}
[c]{ccc}%
\dfrac{\sin(a)}{\sin(\theta)}=\dfrac{\sin(\rho)}{\sin(\beta)}, &  &
\dfrac{\sin(b)}{\sin(\beta)}=\dfrac{\sin(c)}{\sin(\gamma)}%
\end{array}
\]
hence
\begin{align*}
\sin(a)\sin(b)\sin(\gamma)  &  =\dfrac{\sin(\rho)\sin(\theta)}{\sin(\beta
)}\frac{\sin(\beta)\sin(c)}{\sin(\gamma)}\sin(\gamma)\\
\  &  =\sin(\rho)\sin(\theta)\sin(c)\\
\  &  =2\sin(\rho)\sin(\theta)\sin(c/2)\cos(c/2).
\end{align*}
Eriksson \cite{Er} proved that
\[
\tan\left(  \frac\sigma2\right)  =\frac{\sin(a)\sin(b)\sin(\gamma)}%
{1+\cos(a)+\cos(b)+\cos(c)}
\]
from which
\begin{equation}
\tan\left(  \frac\sigma2\right)  =\frac{\sin(c/2)\sin(\rho)\sin(\theta)}%
{\cos(c/2)+\cos(\rho)} \label{eq63}%
\end{equation}
follows immediately.

Adding equation (\ref{eq62}) to (\ref{eq61}), we obtain \cite{Do}
\[
\cos(\rho)=\frac{\cos(a)+\cos(b)}{2\cos\left(  \frac c2\right)  }=\frac
{\cos\left(  \frac{a+b}2\right)  \cos\left(  \frac{a-b}2\right)  }{\cos\left(
\frac c2\right)  };
\]
subtracting equation (\ref{eq62}) from (\ref{eq61}), we obtain
\[
\cos(\theta)=\frac{\cos(a)-\cos(b)}{2\sin\left(  \frac c2\right)  \sin(\rho
)}=-\frac{\sin\left(  \frac{a+b}2\right)  \sin\left(  \frac{a-b}2\right)
}{\sin\left(  \frac c2\right)  \sin(\rho)}.
\]
Thus, given $a,b,c$, it is easy to compute $\rho$ and then $\theta$ (in that order).

Rearranging equation (\ref{eq63}) to
\[
\frac{\tan\left(  \sigma/2\right)  }{\sin(\theta)}=\frac{\sin(c/2)\sin(\rho
)}{\cos(c/2)+\cos(\rho)}=\frac{\sqrt{1-\cos(c/2)^{2}}\sqrt{1-\cos(\rho)^{2}}%
}{\cos(c/2)+\cos(\rho)},
\]
that is,
\[
z=\frac{\sqrt{1-y^{2}}\sqrt{1-x^{2}}}{y+x}
\]
we solve for $x=\cos(\rho)$:
\[
(y+x)^{2}z^{2}=\left(  1-y^{2}\right)  \left(  1-x^{2}\right)  ,
\]
that is
\[
\left(  1-y^{2}+z^{2}\right)  x^{2}+\left(  2yz^{2}\right)  x-\left(
1-y^{2}-y^{2}z^{2}\right)  =0
\]
and obtain the expression for $\rho=f(\theta).$

\subsection{\label{TrigBisector}Angle Bisector Case}

Let $\rho$, $\theta$ be defined within triangle $ABC$ as in \textbf{[}%
\ref{BisectorPerimeter}\textbf{]. }Applying the Law of Cosines for Angles to
both triangles $AQC$ and $AQB$, we have
\begin{equation}
-\cos(\gamma)=\cos(\theta)\cos(\alpha/2)+\sin(\theta)\sin(\alpha/2)\cos(\rho),
\label{eq64}%
\end{equation}
\begin{equation}
-\cos(\beta)=-\cos(\theta)\cos(\alpha/2)+\sin(\theta)\sin(\alpha/2)\cos(\rho)
\label{eq65}%
\end{equation}
because $\cos(\pi-\theta)=-\cos(\theta)$ and $\sin(\pi-\theta)=\sin(\theta)$;
hence
\begin{align*}
\cos(\alpha)+\cos(\beta)+\cos(\gamma)-1  &  =\cos(\alpha)-2\sin(\theta
)\sin(\alpha/2)\cos(\rho)-1\\
\  &  =-2\sin(\theta)\sin(\alpha/2)\cos(\rho)-2\sin(\alpha/2)^{2}\\
\  &  =-2\sin(\alpha/2)\left(  \sin(\theta)\cos(\rho)+\sin(\alpha/2)\right)  .
\end{align*}
Applying the Law of Sines to triangle $AQB$, we have
\[
\frac{\sin(c)}{\sin(\theta)}=\frac{\sin(\rho)}{\sin(\beta)}
\]
hence
\begin{align*}
\sin(\alpha)\sin(\beta)\sin(c)  &  =\sin(\alpha)\frac{\sin(\rho)\sin(\theta
)}{\sin(c)}\sin(c)\\
\  &  =\sin(\alpha)\sin(\rho)\sin(\theta)\\
\  &  =2\sin(\alpha/2)\cos(\alpha/2)\sin(\rho)\sin(\theta).
\end{align*}
The dual of Eriksson's result is \cite{Er}
\[
\tan\left(  \frac\tau2\right)  =\frac{\sin(\alpha)\sin(\beta)\sin(c)}%
{\cos(\alpha)+\cos(\beta)+\cos(\gamma)-1}
\]
from which
\begin{equation}
\tan\left(  \frac\tau2\right)  =-\frac{\cos(\alpha/2)\sin(\rho)\sin(\theta
)}{\sin(\alpha/2)+\cos(\rho)\sin(\theta)} \label{eq66}%
\end{equation}
follows immediately.

Adding equation (\ref{eq65}) to (\ref{eq64}), we obtain \cite{Do}
\begin{align*}
\cos(\rho)  &  =-\frac{\cos(\beta)+\cos(\gamma)}{2\sin\left(  \frac
\alpha2\right)  \sin(\theta)}=-\frac{\cos\left(  \frac{\beta+\gamma}2\right)
\cos\left(  \frac{\beta-\gamma}2\right)  }{\sin\left(  \frac\alpha2\right)
\sin(\theta)};
\end{align*}
subtracting equation (\ref{eq65}) from (\ref{eq64}), we obtain
\[
\cos(\theta)=\frac{\cos(\beta)-\cos(\gamma)}{2\cos\left(  \frac\alpha2\right)
}=-\frac{\sin\left(  \frac{\beta+\gamma}2\right)  \sin\left(  \frac
{\beta-\gamma}2\right)  }{\cos\left(  \frac\alpha2\right)  }.
\]
Thus, given $\alpha,\beta,\gamma$, it is easy to compute $\theta$ and then
$\rho$ (in that order).

Rearranging equation (\ref{eq66}) to
\[
\tan(\tau/2)\sin(\alpha/2)+\left(  \tan(\tau/2)\cos(\rho)+\cos(\alpha
/2)\sin(\rho)\right)  \sin(\theta)=0,
\]
that is,
\[
\sin(\theta)=-\frac{\tan(\tau/2)\sin(\alpha/2)}{\tan(\tau/2)\cos(\rho
)+\cos(\alpha/2)\sin(\rho)}
\]
we obtain the expression for $\theta=f(\rho)$. The smallest admissible value
$\rho=\rho_{\text{thres}}$ occurs when $\sin(\theta)=1$, that is,
\[
\tan(\tau/2)\cos(\rho)+\cos(\alpha/2)\sin(\rho)+\tan(\tau/2)\sin
(\alpha/2)=0\text{.}
\]
Solving
\[
y\,x+z\sqrt{1-x^{2}}+y\sqrt{1-z^{2}}=0
\]
is made possible via
\[
y^{2}\left(  x+\sqrt{1-z^{2}}\right)  ^{2}=z^{2}\left(  1-x^{2}\right)  ,
\]
hence
\[
\left(  y^{2}+z^{2}\right)  x^{2}+2y^{2}\sqrt{1-z^{2}}x+\left(  y^{2}%
-z^{2}-y^{2}z^{2}\right)  =0
\]
hence
\[
x=\frac{-y^{2}\sqrt{1-z^{2}}-z^{2}\sqrt{1+y^{2}}}{y^{2}+z^{2}}=-\frac
{\frac1{\sqrt{1+y^{2}}}+\sqrt{1-z^{2}}}{1+\frac1{\sqrt{1+y^{2}}}\sqrt{1-z^{2}%
}}
\]
hence
\[
\cos(\rho)=-\frac{\cos(\tau/2)+\sin(\alpha/2)}{1+\cos(\tau/2)\sin(\alpha/2)}
\]
gives the desired threshold.

\section{Definite Integrals}

\subsection{\label{CroftonExhum}Crofton/Exhumatus}

We wish to evaluate
\[
\left\{
\begin{array}
[c]{ccc}%
{\displaystyle\int\limits_{0}^{\pi}}
\dfrac{\pi-\arctan\left(  \sqrt{\tan^{2}\frac x2\csc^{2}\frac y2+1}\tan\frac
y2\right)  }{\sqrt{\tan^{2}\frac x2\csc^{2}\frac y2+1}}\sin x\,dx &  &
\text{if }0\leq y<\pi,\\
-%
{\displaystyle\int\limits_{0}^{\pi}}
\dfrac{\arctan\left(  \sqrt{\tan^{2}\frac x2\csc^{2}\frac y2+1}\tan\frac
y2\right)  }{\sqrt{\tan^{2}\frac x2\csc^{2}\frac y2+1}}\sin x\,dx &  &
\text{if }\pi\leq y\leq2\pi.
\end{array}
\right.
\]
A miraculous substitution
\[%
\begin{array}
[c]{ccc}%
\cos z=\cos\tfrac x2\cos\tfrac y2, &  & 0\leq z\leq\pi
\end{array}
\]
is due to Crofton \& Exhumatus \cite{CE}; from
\[
\sin z\,dz=\tfrac12\sin\tfrac x2\cos\tfrac y2\,dx
\]
we deduce that
\[
\sin x\,dx=2\sin\tfrac x2\cos\tfrac x2\frac{\sin z\,dz}{\tfrac12\sin\tfrac
x2\cos\tfrac y2}=4\frac{\cos z}{\cos\tfrac y2}\frac{\sin z\,dz}{\cos\tfrac
y2}=\frac{4\cos z\sin z}{\cos^{2}\tfrac y2}dz
\]
and
\[
\tan^{2}\tfrac x2\csc^{2}\tfrac y2+1=\frac{\tan^{2}z}{\tan^{2}\tfrac y2}
\]
because
\begin{align*}
\tan^{2}z+1  &  =\sec^{2}z=\sec^{2}\tfrac x2\sec^{2}\tfrac y2=\left(  \tan
^{2}\tfrac x2+1\right)  \sec^{2}\tfrac y2\\
&  =\tan^{2}\tfrac x2\sec^{2}\tfrac y2+\left(  \tan^{2}\tfrac y2+1\right) \\
&  =\left(  \tan^{2}\tfrac x2\sec^{2}\tfrac y2+\tan^{2}\tfrac y2\right)  +1\\
&  =\left(  \tan^{2}\tfrac x2\csc^{2}\tfrac y2+1\right)  \tan^{2}\tfrac y2+1.
\end{align*}
Since $\cos(z)$, $\cos(y/2)$ obviously have the same sign, it follows that
$\tan(z)$, $\tan(y/2)$ likewise have the same sign and%

\[
\tan z=\sqrt{\tan^{2}\tfrac x2\csc^{2}\tfrac y2+1}\tan\tfrac y2.
\]

If $0\leq y<\pi$, clearly $\cos(y/2)>0$ and the range $0\leq x\leq\pi$ maps to
$y/2\leq z\leq\pi/2$. Also, $\tan(y/2)>0$, hence
\[
z=\arctan\left(  \sqrt{\tan^{2}\tfrac x2\csc^{2}\tfrac y2+1}\tan\tfrac
y2\right)
\]
hence
\begin{align*}
\dfrac{\pi-\arctan\left(  \sqrt{\tan^{2}\frac x2\csc^{2}\frac y2+1}\tan\frac
y2\right)  }{\sqrt{\tan^{2}\frac x2\csc^{2}\frac y2+1}}\sin x\,dx  &
=\frac{\tan\tfrac y2}{\tan z}(\pi-z)\frac{4\cos z\sin z}{\cos^{2}\tfrac
y2}dz\\
&  =\frac{4\tan\tfrac y2}{\cos^{2}\tfrac y2}(\pi-z)\cos^{2}z\,dz.
\end{align*}
The required definite integral is thus
\[
\frac{4\tan\tfrac y2}{\cos^{2}\tfrac y2}%
{\displaystyle\int\limits_{y/2}^{\pi/2}}
(\pi-z)\cos^{2}z\,dz
\]
which is elementary.

If $\pi<y\leq2\pi$, clearly $\cos(y/2)<0$ and the range $0\leq x\leq\pi$ maps
to $y/2\geq z\geq\pi/2$. Also, $\tan(y/2)<0$, hence
\[
z=\pi+\arctan\left(  \sqrt{\tan^{2}\tfrac x2\csc^{2}\tfrac y2+1}\tan\tfrac
y2\right)
\]
hence
\begin{align*}
-\dfrac{\arctan\left(  \sqrt{\tan^{2}\frac x2\csc^{2}\frac y2+1}\tan\frac
y2\right)  }{\sqrt{\tan^{2}\frac x2\csc^{2}\frac y2+1}}\sin x\,dx  &
=-\frac{\tan\tfrac y2}{\tan z}(z-\pi)\frac{4\cos z\sin z}{\cos^{2}\tfrac
y2}dz\\
\  &  =\frac{4\tan\tfrac y2}{\cos^{2}\tfrac y2}(\pi-z)\cos^{2}z\,dz.
\end{align*}
The required definite integral is thus identical to before (although here the
lower limit $y/2$ is greater than the upper limit $\pi/2$).

Computer algebra swiftly gives
\begin{align*}
&  \ \ \ 1+\frac d{dy}\left[  \frac{4\tan\tfrac y2}{\cos^{2}\tfrac y2}%
{\displaystyle\int\limits_{y/2}^{\pi/2}}
(\pi-z)\cos^{2}z\,dz\right] \\
\  &  =-\frac{(y^{2}-4\pi y+3\pi^{2}-6)\cos(y)-6(y-2\pi)\sin(y)-2(y^{2}-4\pi
y+3\pi^{2}+3)}{8\cos(y/2)^{4}}%
\end{align*}
and this is useful at the conclusion of [\ref{Success}].

\subsection{\label{JonesBenyTink}Jones/Benyon-Tinker}

Combining our results with those in \cite{ABT1, ABT2}, we have
\begin{align}
&  \ \ \ \ \
{\displaystyle\int\limits_{\tau/2-\kappa}^{\tau/2}}
\frac{\sin(\tau-\kappa-\rho)\sin(\kappa)\sin(\rho)}{\sqrt{\sin(\kappa)^{2}%
\sin(\rho)^{2}-\left[  \cos(\kappa)\cos(\rho)-\cos(\tau-\kappa-\rho)\right]
^{2}}}d\rho\label{eq77}\\
\  &  =\frac{E\left(  \sin\left(  \frac\kappa2\right)  \right)  -\cos\left(
\frac{\tau-\kappa}2\right)  ^{2}K\left(  \sin\left(  \frac\kappa2\right)
\right)  }{\sqrt{\cos\left(  \frac\kappa2\right)  ^{2}-\cos\left(  \frac
{\tau-\kappa}2\right)  ^{2}}}\sin(\kappa)\,\nonumber
\end{align}
in connection with primal perimeter [\ref{PrimalPeri}] and
\begin{align}
&  \ \ -\
{\displaystyle\int\limits_{\sigma/2}^{\pi-(\kappa-\sigma/2)}}
\frac{\sin(\sigma-\kappa-\theta)\sin(\kappa)\sin(\theta)}{\sqrt{\sin
(\kappa)^{2}\sin(\theta)^{2}-\left[  \cos(\kappa)\cos(\theta)-\cos
(\sigma-\kappa-\theta)\right]  ^{2}}}d\theta\label{eq78}\\
\  &  =\frac{E\left(  \cos\left(  \frac\kappa2\right)  \right)  -\sin\left(
\frac{\sigma-\kappa}2\right)  ^{2}K\left(  \cos\left(  \frac\kappa2\right)
\right)  }{\sqrt{\sin\left(  \frac\kappa2\right)  ^{2}-\sin\left(
\frac{\sigma-\kappa}2\right)  ^{2}}}\sin(\kappa)\,\nonumber
\end{align}
in connection with dual area [\ref{DualArea}]. A\ direct symbolic proof of
these formulas is not known [\ref{Unsuccess}].

Consider the problem of integrating equation (\ref{eq77}) with respect to
$\kappa$, $0\leq\kappa\leq\tau/2$ and of integrating equation (\ref{eq78})
with respect to $\kappa$, $\sigma/2\leq\kappa\leq\pi$. In (\ref{eq77}), $\rho$
is integrated out first, $\kappa$ second. In (\ref{eq78}), $\theta$ is
integrated out first, $\kappa$ second. By symmetry, we gain nothing by
integrating out $\kappa$ first, thus a closed-form expression for
unconditional density would seem unlikely. Another miraculous change of
variables might, however, be brought into play. Other approaches based on
other coordinate systems exist [\ref{AngleCoord}, \ref{SideCoord}]. It is
still too early to rule out the possibility of a breakthrough here.

\section{Acknowledgement}

We are grateful to M. Larry Glasser for a helpful discussion about integrals
at the end of [\ref{PrimalPeri}] \&\ [\ref{DualArea}]. Much more relevant
material can be found at \cite{F2, F3}, including experimental computer runs
that aided theoretical discussion here. The book \cite{ABT2} studies length
distributions for open and closed random $n$-step tours ($n\geq3$) on spheres,
thus generalizing our discussion of triangle perimeters considerably.

\begin{tabular}
[c]{lllll}
& Steven R. Finch &  & Antonia J. Jones & \\
& Dept. of Statistics &  & School of Computer Science \&\ Informatics & \\
& Harvard University &  & Cardiff University & \\
& Cambridge, MA, USA &  & Cardiff, Wales, UK & \\
& \textit{steven\_finch@harvard.edu} &  & \textit{(deceased)} &
\end{tabular}

\end{document}